# Squircular Calculations

Chamberlain Fong spectralfft@yahoo.com
Joint Mathematics Meetings 2018, SIGMAA-ARTS

**Abstract** – The Fernandez-Guasti squircle is a plane algebraic curve that is an intermediate shape between the circle and the square. It has qualitative features that are similar to the more famous Lamé curve. However, unlike the Lamé curve which has unbounded exponents, the Fernandez-Guasti squircle is a low degree quartic curve. This makes it more amenable to algebraic manipulation and simplification. In this paper, we will analyze this squircle and derive formulas for its area, arc length, and polar form. We will also provide several parametric equations of this squircle. Finally, we extend the Fernandez-Guasti squircle to three dimensions by coming up with an analogous surface that is an intermediate shape between the sphere and the cube.

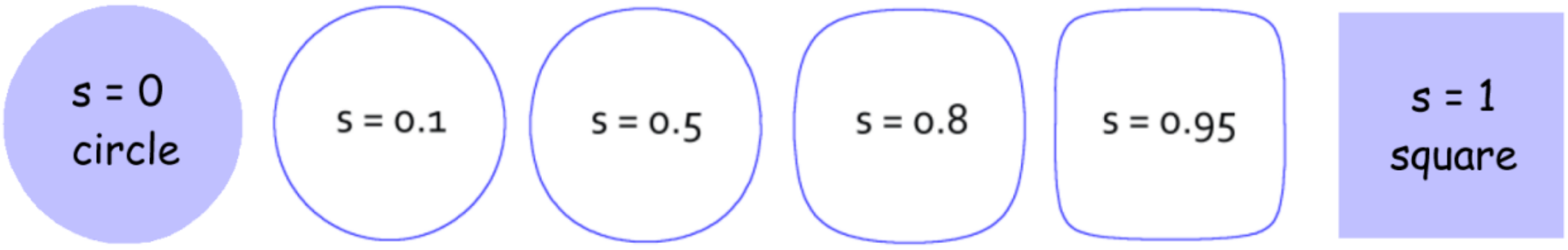

Figure 1: The FG-squircle $x^2 + y^2 - \frac{s^2}{r^2}\, x^2 y^2 = r^2$ at varying values of s.

## 1 Introduction

In 1992, Manuel Fernandez-Guasti discovered an intermediate shape between the circle and the square [Fernandez-Guasti 1992]. This shape is a quartic plane curve. The equation for this shape is provided with Figure 1. There are two parameters for this equation: *s* and *r*. The squareness parameter *s* allows the shape to interpolate between the circle and the square. When *s* = 0, the equation produces a circle with radius *r*. When *s* = 1, the equation produces a square with a side length of *2r*. In between, the equation produces a smooth planar curve that resembles both the circle and the square. For the rest of this paper, we shall refer to this shape as the *FG-squircle*.

In this paper, we shall primarily focus in the region where $-r \le x \le r$ and $-r \le y \le r$. The equation given above can actually have points (x,y) outside our region of interest, but we shall ignore these portions of the curve. For our purpose, the FG-squircle is just a closed curve clipped inside a square region centered at the origin. This clipping region has side length 2r. All points (x,y) outside this clipping square region will be culled out and ignored.

## 2 Derivation of the FG-squircle

Manuel Fernandez-Guasti derived the equation for his squircle by first observing that one could use a product of two binomial square root terms to represent a square shape [Fernandez-Guasti 1992]. Specifically, he noticed that the equation

$$\sqrt{1 - \frac{x^2}{r^2}}\sqrt{1 - \frac{y^2}{r^2}} = 0$$

can be plotted to form a square in the Cartesian plane when $-r \le x \le r$ and $-r \le y \le r$. He then simplified the equation by multiplying out the terms.

$$\sqrt{(1-\frac{x^2}{r^2})(1-\frac{y^2}{r^2})} \;=\; \sqrt{1-\frac{x^2}{r^2}-\frac{y^2}{r^2}+\frac{x^2y^2}{r^4}} \;=0$$

This expression further simplifies to a quartic polynomial equation: $\frac{x^2}{r^2}+\frac{y^2}{r^2}-\frac{x^2y^2}{r^4}=1$

This is a valid polynomial representation of the square in the Cartesian plane when $-r \leq x \leq r$ and $-r \leq y \leq r$. The polynomial is quartic (i.e. degree 4) because of the $\frac{-x^2y^2}{r^4}$ term.

Moreover, Fernandez-Guasti noticed that by introducing a squareness parameter *s* into the equation, one could produce a shape that interpolates between the square and the circle.

$$\frac{x^2}{r^2}+\frac{y^2}{r^2}-\frac{s^2x^2y^2}{r^4}=1$$

$$\Rightarrow \qquad x^2+y^2-\frac{s^2}{r^2}\,x^2y^2=r^2$$

Furthermore, Fernandez-Guasti suggested that eccentricities could be introduced to the equation in order to produce a shape that interpolates between the rectangle and the ellipse. The eccentric version of the equation is

$$\frac{x^2}{a^2}+\frac{y^2}{b^2}-\frac{s^2x^2y^2}{a^2b^2}=1$$

where *a* and *b* in this context are semi-major and semi-minor lengths of the ellipse.

The FG-squircle is closely related to the *Edwards curve* which was studied by Harold Edwards in 2007 [Edwards 2007]. The Edwards curve has the equation

$$x^2+y^2=1+dx^2y^2$$

where *d* is an adjustable scalar parameter. The Edwards curve is equivalent to the FG-squircle of radius *1* when $d=s^2$. In a way, the Edwards curve is a generalized family of curves that encompass the FG-squircle.

On a historical note, Manuel Fernandez-Guasti introduced his squircle in 1992 within the context of plotting curves in Cartesian coordinates for Mathematics education. In contrast, Harold Edwards introduced his curve in 2007 within the context of studying elliptic curves.

More recently, Daniel Bernstein et al. introduced the *twisted Edwards curve* [Bernstein 2008] for use in cryptography. The twisted Edwards curve is an eccentric variant of the Edwards curve with the equation

$$ax^2+y^2=1+d\,x^2y^2$$

## 3 Algebraic Components of the FG-squircle

Although the main focus of this paper with the FG-squircle is inside the region within $-r \leq x \leq r$ and $-r \leq y \leq r$, it is important to mention that the curve has portions outside this region. In fact, the quartic equation for the curve produces five disjoint algebraic components in the 2D Cartesian plane. This is shown in Figure 2. Aside from the central squircular core, there are four infinite components that surround it. We shall call these extraneous components as the *lobes* of the FG-squircle. There is one infinite lobe for each quadrant of the 2D Cartesian plane.

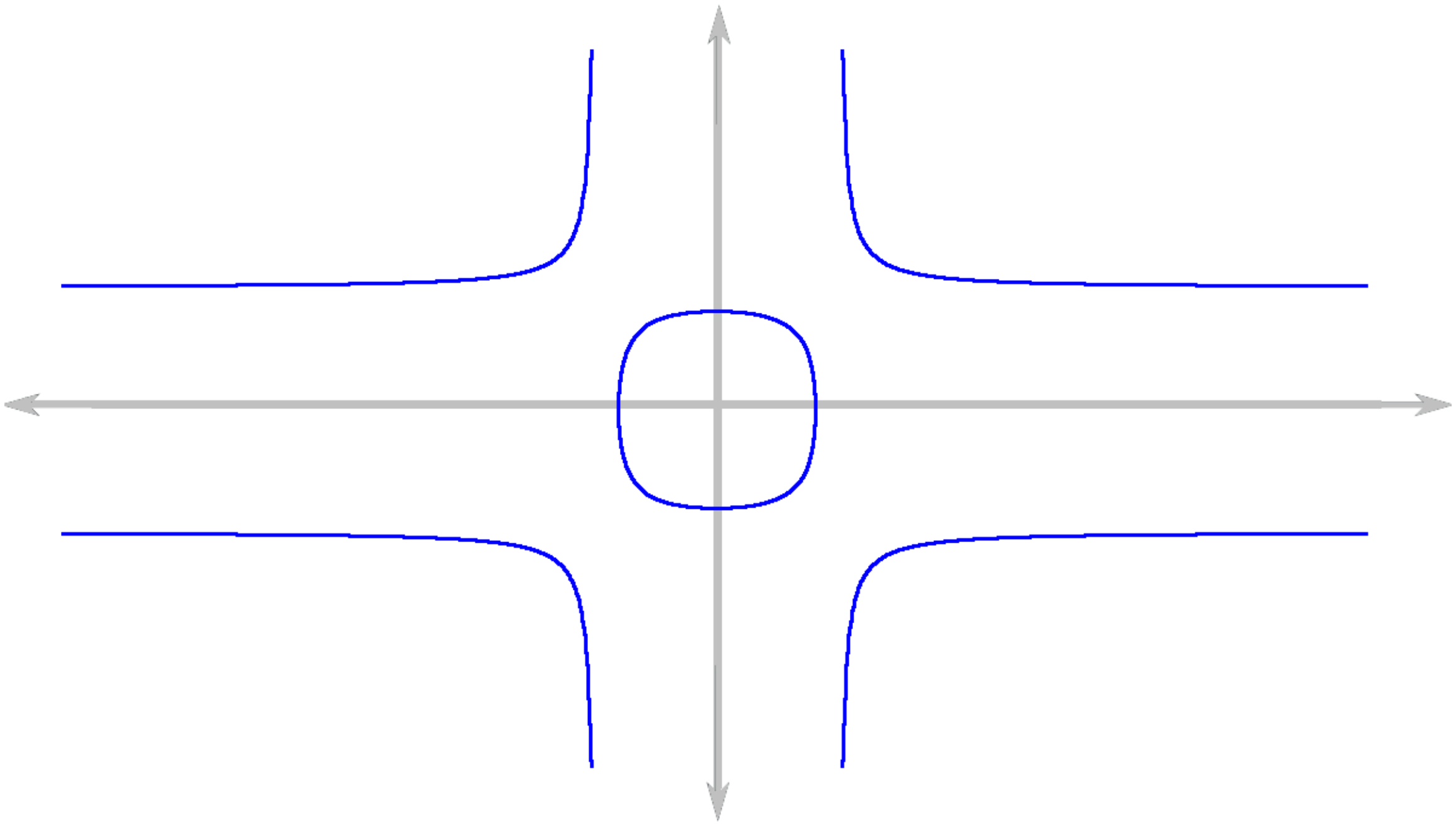

Figure 2: The squircular core curve surrounded by its four lobes
Note: The lobes are disconnected to the central core except when *s=1*

The five disjoint algebraic components of the FG-squircle are contained within separate regions of the 2D Cartesian plane. The bounds are given below:

- Squircular core: $x \in [-r, r]$ and $y \in [-r, r]$
- Lobe#1 (1st quadrant): $x \in \left[\frac{r}{s}, +\infty\right)$ and $y \in \left[\frac{r}{s}, +\infty\right)$
- Lobe#2 (2nd quadrant): $x \in \left(-\infty, \frac{-r}{s}\right]$ and $y \in \left[\frac{r}{s}, +\infty\right)$
- Lobe#3 (3rd quadrant): $x \in \left(-\infty, \frac{-r}{s}\right]$ and $y \in \left(-\infty, \frac{-r}{s}\right]$
- Lobe#4 (4th quadrant: $x \in \left[\frac{r}{s}, +\infty\right)$ and $y \in \left(-\infty, \frac{-r}{s}\right]$

The location and span of the four lobes depend on the squareness parameter of the squircle. When *s=0*, the FG-squircle reduces to a circle and the lobes are effectively infinitely far away; i.e. gone. For *s >0*, each of the four lobes will appear at a finite distance from the squircular core. As *s* increases, the four lobes will get closer and closer to the squircular core. When *s=1*, the four lobes will touch the squircular core. This will result in four points of self-intersection at the corners of the squircle. When *s=1*, the five disjoint algebraic components of the squircle will join into one connected curve with four singularities at the corners of the square.

## 4 The Lamé Curve

In 1818, Gabriel Lamé studied a plane curve that now bears his name. Qualitatively, it appears very similar to the FG-squircle, but they have important differences which we will discuss later. The Lamé curve has the equation

$$|x|^P + |y|^P = r^P$$

There are two parameters for this equation: *p* and *r*. The power parameter *p* allows the shape to interpolate between the circle and the square. The parameter *r* determines the size of the resulting shape. The Lamé curve actually contains two families of squircles -- the Lamé upper squircle and the Lamé lower squircle. These are shown in Figure 3.

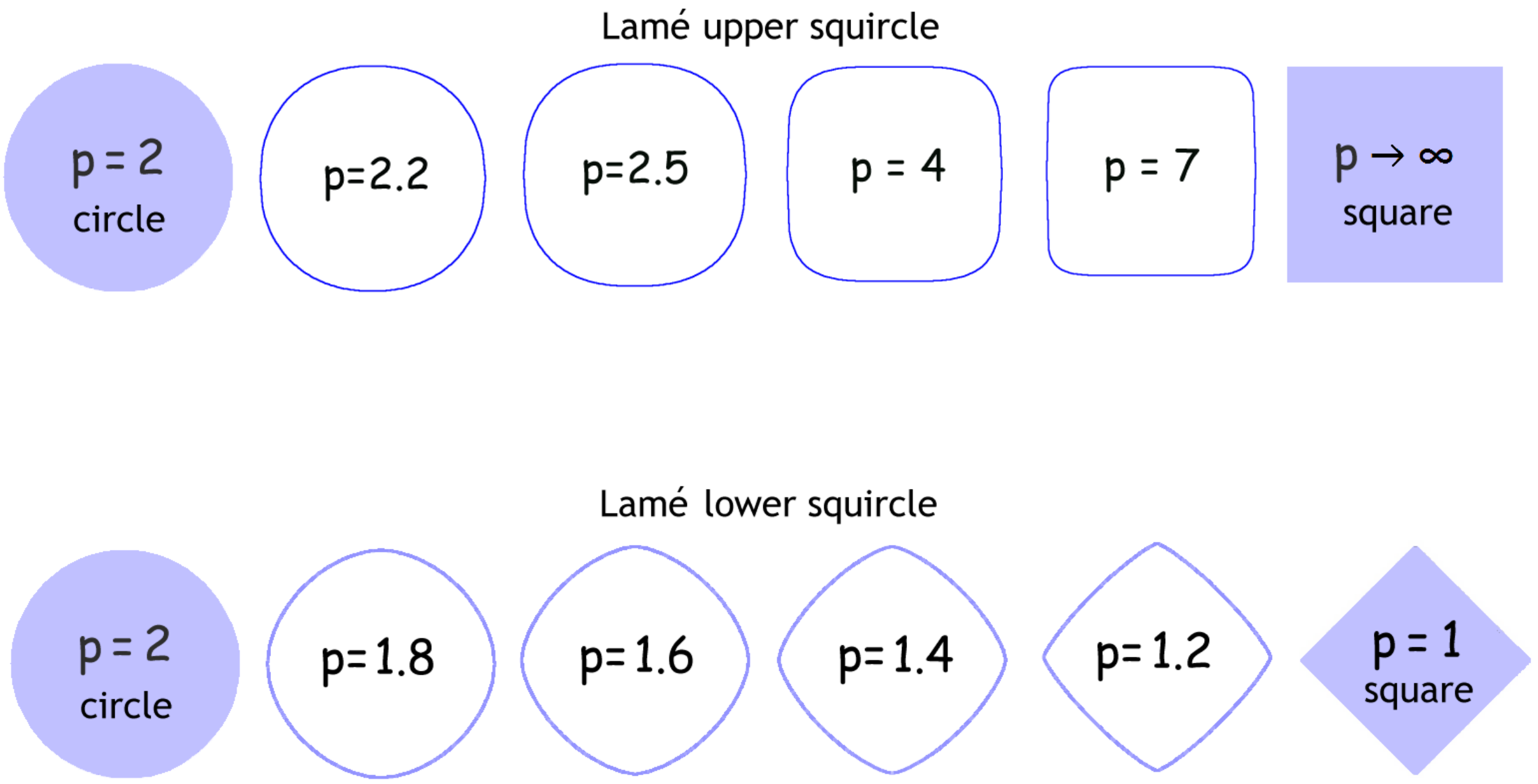

Figure 3: The Lamé curve at varying values of p.

The *Lamé upper squircle* results from the equation for $p \in [2,+\infty)$. When $p = 2$, the equation produces a circle with radius *r*. As $p \to \infty$, the equation produces a square with a side length of *2r*. In between, the equation produces a smooth planar curve that resembles both the circle and the square. The upper portion of Figure 3 shows the Lamé upper squircle at varying powers.

The *Lamé lower squircle* results from the equation for $p \in [1,2]$. When $p = 2$, the equation produces a circle with radius *r*. When $p = 1$, the equation produces a tilted square with a side length of $r\sqrt{2}$. In between, the equation produces a smooth planar curve that resembles both the circle and the square. The lower portion of Figure 3 shows the Lamé lower squircle at varying powers.

Although the Lamé curve exhibits qualitative similarities with the FG-squircle, there is one major difference between the two. The Lamé curve can only approximate an upright square. It requires an infinite exponent in order to fully realize an upright square. Moreover, since the equation for the Lamé curve has unbounded exponents, it is unwieldy and difficult to manipulate algebraically.

In contrast, the FG-squircle is a low-degree quartic curve. In addition, since there are no cubic or linear terms in its polynomial equation, the FG-squircle is essentially an even simpler biquadratic curve. This makes it much more amenable to algebraic manipulation and simplification. In fact, because of its algebraic simplicity, the FG-squircle has been used to produce explicit and invertible mappings [Fong 2014][Fong 2017] between the circular disc and the square. Furthermore, in sections 10 and 15 of this paper, we shall show 3D algebraic surfaces with simple polynomial equations that are based on the FG-squircle. This is possible because of the algebraic simplicity of the FG-squircle.

Meanwhile, it is now appropriate to clarify our terminology. Currently, there is no standard definition for the word 'squircle'. There is mathematical literature that refer to the squircle as the Lamé curve for $p = 4$; i.e.

$$x^4 + y^4 = r^4$$

We prefer to call to this specific curve as *Lamé's special quartic*. For this paper, the word 'squircle' is an umbrella term for an intermediate shape between the circle and the square. So far, we have discussed three different types of squircle:

- Fernandez-Guasti squircle
- Lamé upper squircle
- Lamé lower squircle

To add to the confusion, the Lamé curve is also referred to as the *superellipse*. However, we prefer to consider the superellipse as the eccentric version of the Lamé curve, with the equation

$$\left|\frac{x}{a}\right|^P + \left|\frac{y}{b}\right|^P = 1$$

## 5 Linearizing the Blending Parameter

Note that using the squareness parameter *s* for blending between the circle and the square produces a non-linear interpolation of the two shapes. One can see in Figure 1 that when s=0.5, the FG-squircle still resembles the circle with very little hint of the square. In fact, the FG-squircle does not resemble a midway shape between the circle and the square until $s \approx 0.928$

It is better to introduce a parameter $\tau \in [0,1]$ that produces a linearized blending of the circle with the square. The simplest way to do this is to let this parameter $\tau$ (Greek letter tau) linearly interpolate the top right corner points of circle and the square.

The top right corner point of the circle is $(\frac{r\sqrt{2}}{2}, \frac{r\sqrt{2}}{2})$. In contrast, the top right corner point of its circumscribing square is *(r,r)* . To produce a diagonal line segment that linearly interpolates between these two corner points with parameter $\tau$, we have these equations:

$$x = \frac{r\sqrt{2}}{2}\,(1-\tau) + \tau\, r$$

$$y = \frac{r\sqrt{2}}{2}\,(1-\tau) + \tau\, r$$

Now if we isolate *s* in the FG-squircle equation, we get

$$x^2 + y^2 - \frac{s^2}{r^2}\, x^2y^2 = r^2 \qquad \Longrightarrow \qquad \frac{s^2x^2y^2}{r^2} = \; x^2 + y^2 - r^2 \qquad \Longrightarrow \qquad s^2 = \frac{r^2}{x^2y^2}\,(x^2 + y^2 - r^2)$$

which simplifies to

$$s = \frac{r}{xy}\sqrt{x^2 + y^2 - r^2}$$

We can then put our corner point linear interpolation equations for *(x,y)* into the expression for *s* to get

$$s = \frac{r}{r^2\left(\frac{\sqrt{2}}{2}\,(1-\tau)+\tau\right)^2}\sqrt{2r^2\left(\frac{\sqrt{2}}{2}\,(1-\tau)+\tau\right)^2 - r^2}$$

Using some algebraic manipulation and simplifying, we get this expression for *s* in terms of $\tau$

$$s = 2\,\frac{\sqrt{(3-2\sqrt{2})\tau^2 - (2-2\sqrt{2})\tau}}{\left(1-(1-\sqrt{2})\tau\right)^2}$$

This gives us a linear blending parameter between the circle and the square. We can then rewrite the equation for the FG-squircle as

$$x^2 + y^2 - \;4\frac{x^2y^2}{r^2}\;\frac{(3-2\sqrt{2})\tau^2 - (2-2\sqrt{2})\tau}{\left(1-(1-\sqrt{2})\tau\right)^4} = r^2$$

This equation of the FG-squircle has parameter $\tau$ instead of *s*. It has the property that as parameter $\tau$ goes linearly from 0 to 1, the top right corner point of the FG-squircle goes linearly from $(\frac{r\sqrt{2}}{2},\frac{r\sqrt{2}}{2})$ to *(r,r)*.

Using this form of the FG-squircle, we can now find a shape midway between circle and the square. This occurs at the value $\tau = 0.5$ . This midway squircle is shown in Figure 4.

We can solve for the midway squareness value as

$$\boldsymbol{s}_{midway} = 2\,\frac{\sqrt{(3-2\sqrt{2})0.5^2-(2-2\sqrt{2})0.5}}{\left(1-(1-\sqrt{2})0.5\right)^2} = \sqrt{736\sqrt{2}-1040} \quad \approx 0.927998872$$

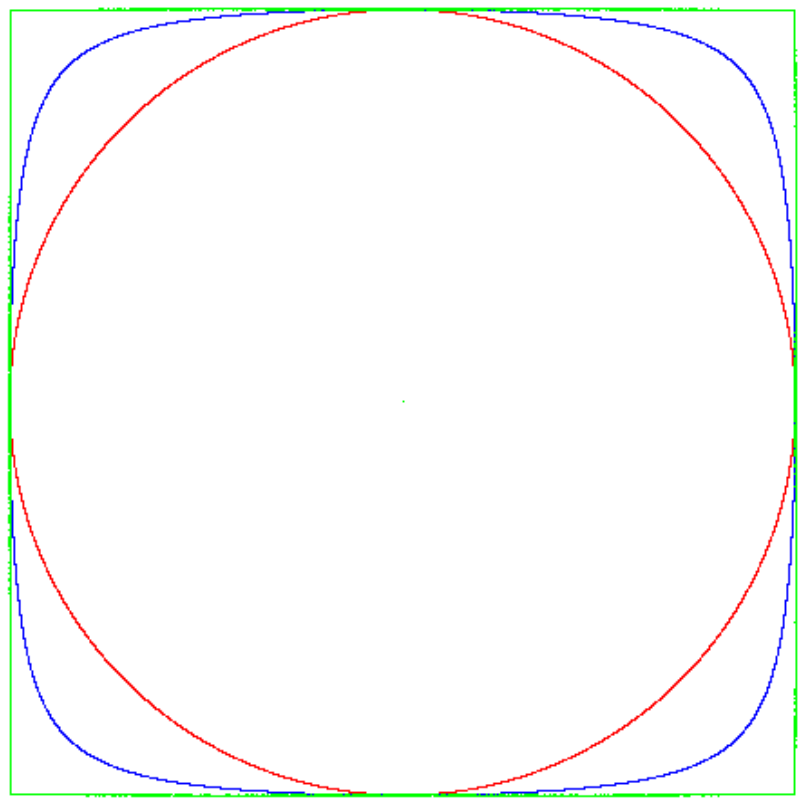

Figure 4: The midway squircle (blue) between the circle (red) and the square (green) at $s \approx 0.928$

## 6 Area of the FG-squircle

We shall derive the formula for the area of the FG-squircle here. Before we proceed with our derivation, we would like to provide a cursory review of the Legendre elliptic integrals.

### 6.1 Legendre Elliptic Integral of the 1st kind

According to chapter 19.2 of the "Digital Library of Mathematical Functions" published by the National Institute of Standards and Technology, the definition of the incomplete Legendre elliptic integral of the 1st kind is:

$$F(\phi, k) = \int_0^{\phi} \frac{d\theta}{\sqrt{1 - k^2 \sin\theta}} = \int_0^{\sin\phi} \frac{dt}{\sqrt{1 - t^2}\sqrt{1 - k^2 t^2}}$$

The argument ϕ is known as the amplitude; and the argument k is known as the modulus. In its standard form, the modulus has a value between 0 and 1, i.e. $0 \le k \le 1$

The complete Legendre elliptic integral of the 1st kind, *K(k)* is defined from its incomplete form with an amplitude $\frac{\pi}{2}$, i.e. $K(k) = F\left(\frac{\pi}{2}, k\right)$

### 6.2 Legendre Elliptic Integral of the 2nd kind

According to chapter 19.2 of the "Digital Library of Mathematical Functions" published by the National Institute of Standards and Technology, the definition of the incomplete Legendre elliptic integral of the 2nd kind is:

$$E(\phi, k) = \int_0^{\phi} \sqrt{1 - k^2 \sin\theta}\, d\theta = \int_0^{\sin\phi} \frac{\sqrt{1 - k^2 t^2}}{\sqrt{1 - t^2}} dt$$

If x = sin ϕ, then we have this formula:

$$E(\sin^{-1} x, k) = \int_0^{x} \frac{\sqrt{1 - k^2 t^2}}{\sqrt{1 - t^2}} dt$$

The complete Legendre elliptic integral of the 2nd kind *E(k)* can be defined from its incomplete form with an amplitude as $\frac{\pi}{2}$ , i.e. $E(k) = E\left(\frac{\pi}{2}, k\right)$

### 6.3 Reciprocal Modulus

As part of our derivation of the area of the FG-squircle, we will need the formula for the reciprocal of the modulus for the Legendre elliptic integral of the 2nd kind. This is provided in chapter 19.7 of the "Digital Library of Mathematical Functions" as

$$E\left(\phi, \frac{1}{q}\right) = \frac{1}{q}\left[E(\beta, q) - (1 - q^2) F(\beta, q)\right]$$

where $\sin\beta = \frac{1}{q}\sin\phi$ or equivalently $\beta = \sin^{-1}\frac{\sin\phi}{q}$

### 6.4 Derivation of Area (incomplete type)

Starting with the equation for the FG-squircle and rearranging the equation to isolate y on one side of the equation, we get:

$$x^2 + y^2 - \frac{s^2}{r^2} x^2 y^2 = r^2 \quad \Rightarrow \quad y^2\left(1 - \frac{s^2}{r^2} x^2\right) = r^2 - x^2 \quad \Rightarrow \quad y = \sqrt{\frac{r^2 - x^2}{1 - \frac{s^2}{r^2} x^2}} = r\sqrt{\frac{1 - \frac{1}{r^2} x^2}{1 - \frac{s^2}{r^2} x^2}}$$

Since the FG-squircle is a shape with dihedral symmetry, the area of the whole FG-squircle is just four times its area in one quadrant, i.e.

$$A = 4\int_0^r y\,dx = 4\int_0^r r\sqrt{\frac{1-\frac{1}{r^2}x^2}{1-\frac{s^2}{r^2}x^2}}\,dx$$

Now using a substitution of variables: $t^2 = \frac{s^2}{r^2}x^2 \quad \Rightarrow \quad dt = \frac{s}{r}\,dx \quad \Rightarrow \quad dx = \frac{r}{s}\,dt$

$$A = 4r\int_0^s \frac{r}{s}\sqrt{\frac{1-\frac{1}{s^2}t^2}{1-t^2}}\,dt = 4\frac{r^2}{s}\int_0^s \frac{\sqrt{1-\frac{1}{s^2}t^2}}{\sqrt{1-t^2}}\,dt$$

Using the equation for the Legendre elliptic integral of the 2nd kind in section 6.2, we can write the area as

$$A = 4\frac{r^2}{s}E(\sin^{-1}s, \tfrac{1}{s})$$

### 6.5 Derivation of Area (complete type)

The areal equation provided in section 6.4 is valid. In fact, this is the formula provided in the Wolfram Mathworld website [Weisstein 2016b]. However, since the squareness parameter $s$ has a value between 0 and 1, its reciprocal will have a value greater than 1. Hence the areal formula from section 6.4 has a modulus in non-standard form. Many software implementations of this Legendre elliptical integral will fail when the modulus k>1. In this section, we will provide an alternative formula for the area of the FG-squircle that does not have this limitation.

Using the reciprocal modulus formula in section 6.3 and substituting q = s = sin ϕ, we get the formula for area as

$$A = 4\frac{r^2}{s}E(\sin^{-1}s, \tfrac{1}{s}) = 4\frac{r^2}{s^2}[E(\sin^{-1}1, s) - (1-s^2)F(\sin^{-1}1, s)] = 4\frac{r^2}{s^2}[E(\tfrac{\pi}{2}, s) + (s^2-1)F(\tfrac{\pi}{2}, s)]$$

This further simplifies by using the definitions of the complete Legendre elliptic integrals.

$$A = 4\frac{r^2}{s^2}\,[\,E(s) + (s^2-1)K(s)]$$

The correctness of this formula for area can be verified numerically by using Monte Carlo methods [Cheney 1999].

## 7 Polar Form of the FG-squircle

In this section, we will derive the equation for the FG-squircle in polar coordinates. In other words, we want to write its equation in the form of $\rho = f(\theta)$ where $\rho = \sqrt{x^2+y^2}$ and $\tan\theta = \frac{y}{x}$

Starting with standard equation for the FG-squircle, $x^2 + y^2 - \frac{s^2}{r^2}x^2y^2 = r^2$

Substitute x and y with their polar form

$$x = \rho\cos\theta$$
$$y = \rho\sin\theta$$

to get

$$\rho^2\cos^2\theta + \rho^2\sin^2\theta - \frac{s^2}{r^2}\rho^4\cos^2\theta\,\sin^2\theta = r^2$$

$$\Rightarrow \quad \rho^2 - \frac{s^2}{r^2}\rho^4\cos^2\theta\,\sin^2\theta = r^2$$

$$\Rightarrow \qquad \rho^2 - \frac{s^2\rho^4}{4r^2}\sin^2 2\theta = r^2$$

$$\Rightarrow \qquad \rho^4 \frac{s^2}{4r^2}\sin^2 2\theta - \rho^2 + r^2 = 0$$

Now, we need to find an equation for $\rho$ in terms of θ. The equation above is a quartic polynomial equation in $\rho$ with no cubic or linear terms. This type of equation is known as a biquadratic and could be solved using the standard quadratic equation.

$$a = \frac{s^2}{4r^2}\sin^2 2\theta \qquad b = -1 \qquad c = r^2$$

$$\Rightarrow \qquad \rho^2 = \frac{1 - \sqrt{1 - s^2\sin^2 2\theta}}{2\frac{s^2}{4r^2}\sin^2 2\theta}$$

Taking the square root of this gives a polar equation of the FG-squircle

$$\boldsymbol{\rho = \frac{r\sqrt{2}}{s\,|\sin 2\theta|}\sqrt{1 - \sqrt{1 - s^2\sin^2 2\theta}}}$$

Manuel Fernandez-Guasti actually wrote this polar equation in his original squircle paper [Fernandez-Guasti 1992]. However, he focused solely on the special case for the square when *s=1*, i.e.

$$\rho = \frac{r\sqrt{2}}{|\sin 2\theta|}\sqrt{1 - \sqrt{1 - \sin^2 2\theta}}$$

Note that this form of the polar equation produces the indeterminate form $\frac{0}{0}$ when $\theta = \frac{n\pi}{2}$, for any integer *n*. However, the limiting value of the equation does exist for these $\theta$, i.e.

$$\lim_{\theta \to \frac{1}{2}n\pi} \frac{r\sqrt{2}}{s\,|\sin 2\theta|}\sqrt{1 - \sqrt{1 - s^2\sin^2 2\theta}} = r$$

One way to mitigate this indeterminate form in software applications such as plotting the FG-squircle in polar coordinates is to introduce a small epsilon value such that $\varepsilon > 0$ to the polar equation

$$\rho = \frac{r\sqrt{2}}{s\,|\sin(2\theta + \varepsilon)|}\sqrt{1 - \sqrt{1 - s^2\sin^2(2\theta + \varepsilon)}}$$

For example, setting $\varepsilon = 0.001$ works well in practice.

However, this polar equation can actually be simplified further using a conjugation trick! Consider multiplying the polar equation with this conjugated fraction

$$\frac{\sqrt{1 + \sqrt{1 - s^2\sin^2 2\theta}}}{\sqrt{1 + \sqrt{1 - s^2\sin^2 2\theta}}}$$

This gives us

$$\rho = \frac{r\sqrt{2}}{s\,|\sin 2\theta|}\sqrt{1 - \sqrt{1 - s^2\sin^2 2\theta}}\;\frac{\sqrt{1 + \sqrt{1 - s^2\sin^2 2\theta}}}{\sqrt{1 + \sqrt{1 - s^2\sin^2 2\theta}}}$$

$$\Rightarrow \qquad \rho = \frac{r\sqrt{2}}{s\,|\sin 2\theta|} \frac{\sqrt{1 - \left(\sqrt{1 - s^2 \sin^2 2\theta}\right)^2}}{\sqrt{1 + \sqrt{1 - s^2 \sin^2 2\theta}}}$$

$$\Rightarrow \qquad \rho = \frac{r\sqrt{2}}{s\,|\sin 2\theta|} \frac{\sqrt{1 - 1 + s^2 \sin^2 2\theta}}{\sqrt{1 + \sqrt{1 - s^2 \sin^2 2\theta}}}$$

$$\Rightarrow \qquad \rho = \frac{r\sqrt{2}}{s\,|\sin 2\theta|} \frac{\sqrt{s^2 \sin^2 2\theta}}{\sqrt{1 + \sqrt{1 - s^2 \sin^2 2\theta}}}$$

$$\Rightarrow \qquad \boldsymbol{\rho = \frac{r\sqrt{2}}{\sqrt{1 + \sqrt{1 - s^2 \sin^2 2\theta}}}}$$

This gives us a much simpler polar equation for the FG-squircle. Moreover, the radical in the numerator can be removed by multiplying the equation with

$$\frac{\sqrt{2}}{\sqrt{2}}$$

to get

$$\rho = \frac{2r}{\sqrt{2 + 2\sqrt{1 - s^2 \sin^2 2\theta}}}$$

which is an alternative form of the polar equation. This simplified form of the polar equation is preferable because it gets rid of the indeterminate form $\frac{0}{0}$ when $\theta = \frac{n\pi}{2}$

Furthermore, a special case of this equation when *s=1* gives us a simple polar equation for the square:

$$\rho = \frac{r\sqrt{2}}{\sqrt{1 + \sqrt{1 - \sin^2 2\theta}}}$$

## 8 Parametric Equations for the FG-squircle

In this section, we will present 6 sets of parametric equations for the FG-squircle. These equations are particularly useful for plotting the FG-squircle. These equations can also be used in the numerical calculation of the arc length of the FG-squircle, which we will discuss in Section 9.2.

### 8.1 Closed Parametric Equations

In 2014, Fong used a simplified FG-squircle to come up with two different methods to map the closed circular disc to a closed square [Fong 2014]. The two mappings are named the *FG-squircular mapping* and the *elliptical grid mapping*, respectively and shown in Figure 5. By *closed*, we mean the shape region including its boundary curve. The forward and inverse equations of the mappings are included in the center column of the figure. In order to illustrate the effects of the mapping, a circular disc with a radial grid is converted to a square and shown in the left side of the figure. Similarly, a square region with a rectangular grid is converted to a circular disc and shown in the right side of the figure.

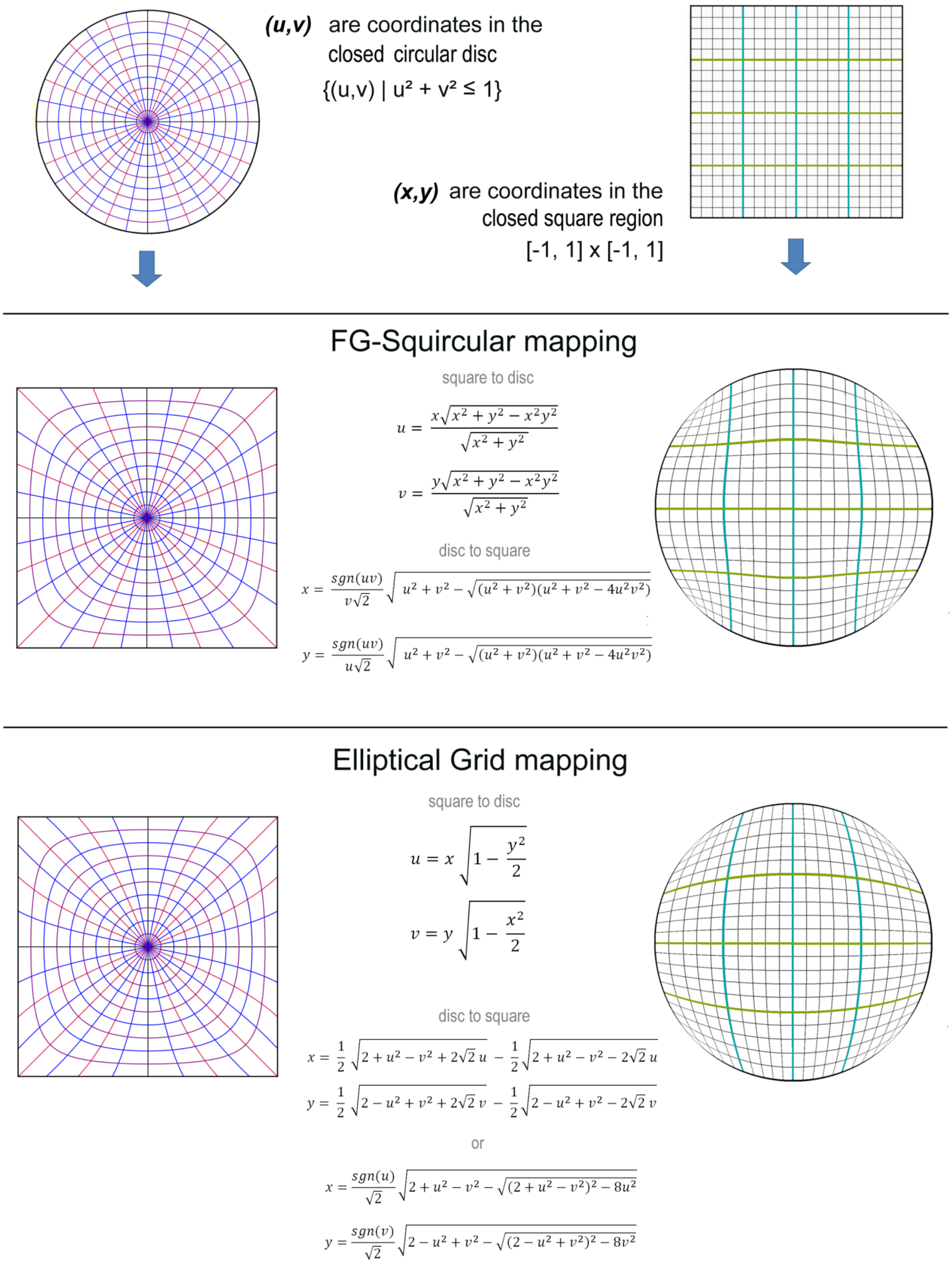

Figure 5: Closed disc-to-square mappings based on the FG-Squircle

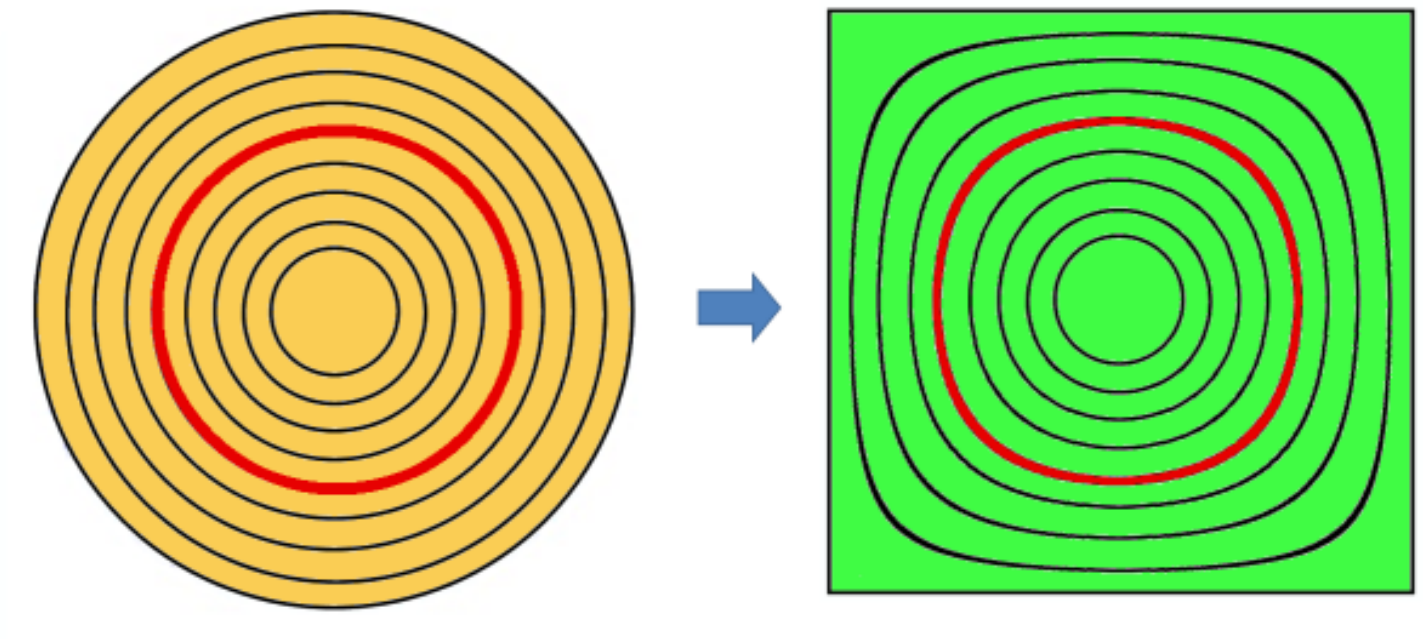

Figure 6: Continuum of circles inside the disc (left)
and continuum of squircles inside the square (right)

Both of these mappings convert concentric circular contours inside the disc to concentric squircular contours inside the square. This property is illustrated in Figure 6. Using this property and the disc-to-square equations in Figure 5, we can come up with parametric equations for the FG-squircle. Just substitute *u=cos(t)* and *v=sin(t)* in the equations then simplify in order to get equations for *(x,y)* in terms of parameter t. Here are the 3 sets of parametric equations for the FG-squircle:

1) Parametric equation based on the *elliptic grid mapping*

$$x = \frac{r}{2s}\sqrt{2 + 2s\sqrt{2}\ \cos t + s^2 \cos 2t} - \frac{r}{2s}\sqrt{2 - 2s\sqrt{2}\ \cos t + s^2 \cos 2t}$$

$$y = \frac{r}{2s}\sqrt{2 + 2s\sqrt{2}\ \sin t - s^2 \cos 2t} - \frac{r}{2s}\sqrt{2 - 2s\sqrt{2}\ \sin t - s^2 \cos 2t}$$

2) Parametric equation based on the *elliptic grid mapping* (alternative)

$$x = \frac{r\ sgn(\cos t)}{s\sqrt{2}}\sqrt{2 + s^2 \cos 2t - \sqrt{(2 + s^2 \cos 2t)^2 - 8\, s^2 \cos^2 t}}$$

$$y = \frac{r\ sgn(\sin t)}{s\sqrt{2}}\sqrt{2 - s^2 \cos 2t - \sqrt{(2 - s^2 \cos 2t)^2 - 8\, s^2 \sin^2 t}}$$

3) Parametric equation based on the *FG-squircular mapping*

$$x = \frac{r\ sgn(\cos t)}{s\sqrt{2}\,|\sin t|}\sqrt{1 - \sqrt{1 - s^2 \sin^2 2t}}$$

$$y = \frac{r\ sgn(\sin t)}{s\sqrt{2}\,|\cos t|}\sqrt{1 - \sqrt{1 - s^2 \sin^2 2t}}$$

Note that parametric equation #3 is equivalent to the polar form of the FG-squircle. This can be verified by using the expression for $\rho$ in Section 7 and substituting into the basic equations for polar coordinates

$$x = \rho \cos \theta \qquad and \qquad y = \rho \sin \theta$$

to get

$$x = \frac{r\sqrt{2}\ \cos t}{s\ \ \sin 2t}\sqrt{1 - \sqrt{1 - s^2 \sin^2 2t}} \qquad and \qquad y = \frac{r\sqrt{2}\sin t}{s\ \ \sin 2t}\sqrt{1 - \sqrt{1 - s^2 \sin^2 2t}}$$

and simplifying by using the trigonometric identity: $sin\ 2t = 2\ sin\ t\ cos\ t$

## 8.2 Open Parametric Equations

In 2017, Fong used the FG-squircle to come up with more methods to map the circular disc to a square. This time, the mappings are *open*, meaning that the mappings do not include the boundary of the circle or square; i.e., the mapping equations only hold for points in the interiors of the circular disc and the square.

The open circular disc is given as the set $\{(u,v) |\ u^2 + v^2 < 1\}$. Note that there is a strict inequality in the expression, which means that the boundary circle is not included. Similarly, the open square region is given as the set *(-1,1) x (-1,1)*. This region also does not include its boundary square.

The mappings are named *squelched grid mapping*, *vertical squelch mapping*, and *horizontal squelch mapping*. These mappings are shown in Figure 7 together with corresponding forward and inverse equations.

Just as we did with the closed mappings in Figure 5, we can come up with parametric equations for the FG-squircle by substituting *u=cos(t)* and *v=sin(t)* in the mapping equations given in Figure 7.

4) Parametric equation based on the *vertical squelch mapping*

$$x = r\cos t$$

$$y = \frac{\mathrm{r}\sin t}{\sqrt{1 - s^2\cos^2 t}}$$

5) Parametric equation based on the *horizontal squelch mapping*

$$x = \frac{r\cos t}{\sqrt{1 - s^2\sin^2 t}}$$

$$y = \mathrm{r}\sin t$$

6) Parametric equation based on the *squelched grid mapping*

$$x = \frac{r\cos t}{\sqrt{1 - \mathrm{s}_1\sin^2 t}}$$

$$y = \frac{\mathrm{r}\sin t}{\sqrt{1 - \mathrm{s}_1\cos^2 t}}$$

$$where \;\; s_1 = s\sqrt{2 - s^2}$$

Here, $s_1$ is a kind of squareness parameter but it is not identical as the *s* squareness parameter in the Fernandez-Guasti squircle equation. The square parameter *s* has this relationship with $s_1$

$$s = \sqrt{1 - \sqrt{1 - s_1^2}}$$

Note that the parametric equations #4 and #5 are asymmetrical; i.e. the expressions for *x* and *y* do not share a similar form that are interchangeable in terms of *cos(t)* and *sin(t)* . In contrast, the parametric equations #1 #2 #3 and #6 are symmetrical with respect to *cos(t)* and *sin(t).*

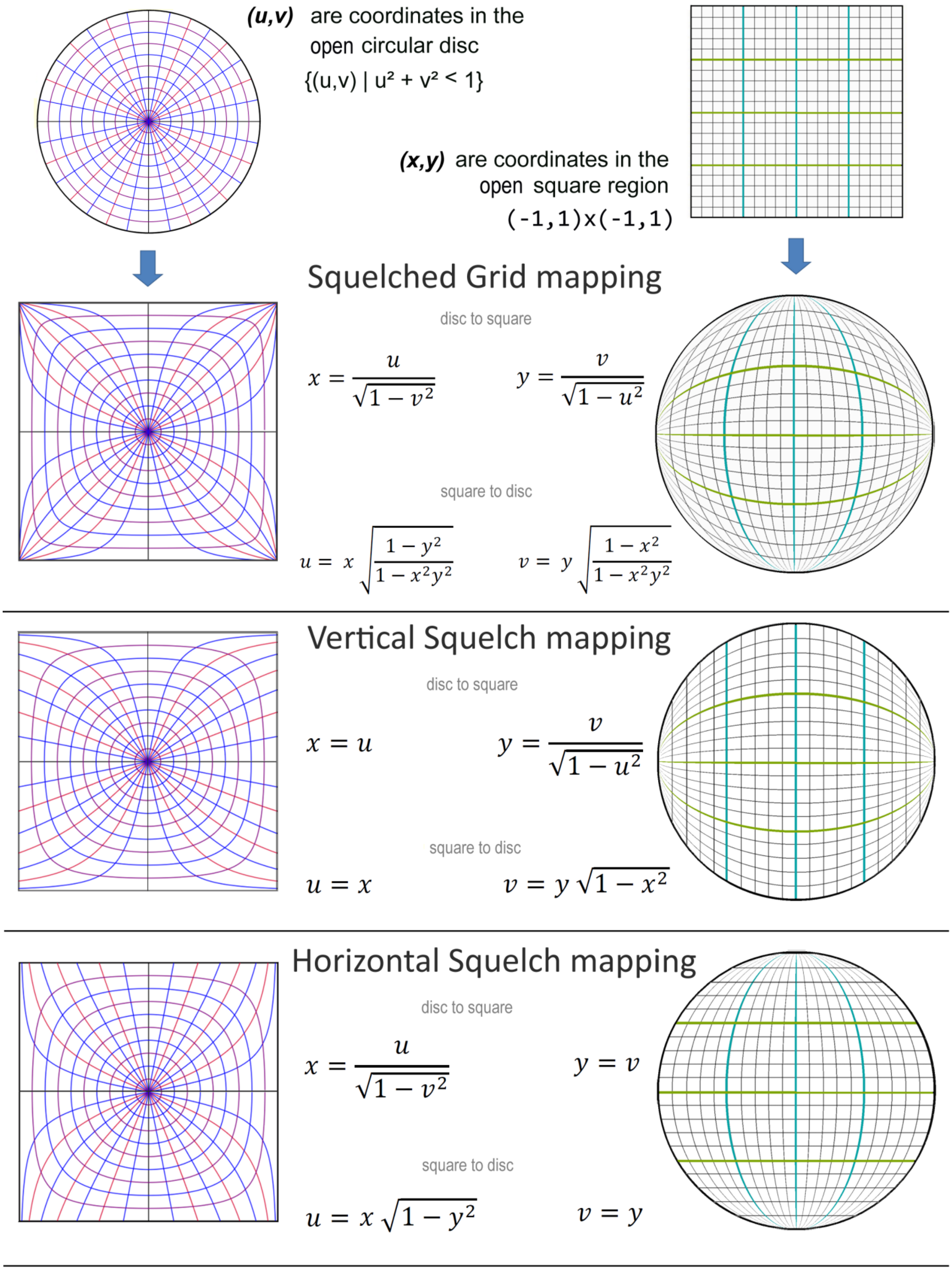

Figure 7: Open disc-to-square mappings based on the FG-Squircle

### 8.3 Parametric Equations via Stereographic Projection

It is known that one can get rational parametric equations for the circle by using the stereographic projection in geometry. For a circle of radius *r* centered at the origin, these parametric equations are

$$x = r\frac{1 - t^2}{1 + t^2}$$

$$y = r\frac{2t}{1 + t^2}$$

Here the parameter t goes from negative infinity to positive infinity, i.e. $t \in (-\infty, +\infty)$. These equations provide a rational formula for the parameterization of the circle that does not involve any trigonometric functions. Using the vertical squelch and horizontal squelch mappings shown in Figure 7, one can get two more sets of parametric equations for the FG-squircle. These parametric equations are also based on the stereographic projection.

Using the vertical squelch mapping, the parametric equations for the FG-squircle becomes

$$x = r\frac{1 - t^2}{1 + t^2} \qquad y = r\frac{2t}{(1 + t^2)\sqrt{1 - s^2x^2}}$$

If we plug-in the *t* formula for *x* into the equation for *y*, the parametric equations further simplify to

$$\boldsymbol{x = r\frac{1 - t^2}{1 + t^2}}$$

$$\boldsymbol{y = r\frac{2t}{\sqrt{(1 + t^2)^2 - s^2(1 - t^2)^2}}}$$

Similarly, using the horizontal squelch mapping, the parametric equations for the FG-squircle becomes

$$x = r\frac{1 - t^2}{(1 + t^2)\sqrt{1 - s^2y^2}} \qquad y = r\frac{2t}{1 + t^2}$$

If we plug-in the *t* formula for *y* into the equation for *x*, the parametric equations simplify to

$$\boldsymbol{x = r\frac{1 - t^2}{\sqrt{(1 + t^2)^2 - 4s^2t^2}}}$$

$$\boldsymbol{y = r\frac{2t}{1 + t^2}}$$

For both sets of parametric equations for the FG-squircle, the parameter t goes from negative infinity to positive infinity, i.e. $t \in (-\infty, +\infty)$.

### 8.4 Parametric Equations for the Square

Using these parametric equations for the Fernandez Guasti squircle, we can easily get parametric equations for the square by setting *s=1*. This square will have a side length of *2r*.

1) Square parametric equation from Elliptic Grid mapping

$$x = \frac{r}{2}\sqrt{2 + 2\sqrt{2}\ \cos t + \cos 2t} - \frac{r}{2}\sqrt{2 - 2\sqrt{2}\ \cos t + \cos 2t}$$

$$y = \frac{r}{2}\sqrt{2 + 2\sqrt{2}\ \sin t - \cos 2t} - \frac{r}{2}\sqrt{2 - 2\sqrt{2}\ \sin t - \cos 2t}$$

2) Square parametric equation from Elliptic Grid mapping (alternative)

$$x = \frac{r\ \ sgn(\cos t)}{\sqrt{2}}\sqrt{2 + \cos 2t - \sqrt{(2 + \cos 2t)^2 - 8\cos^2 t}}$$
$$y = \frac{r\ \ sgn(\sin t)}{\sqrt{2}}\sqrt{2 - \cos 2t - \sqrt{(2 - \cos 2t)^2 - 8\ \sin^2 t}}$$

3) Square parametric equation from FG-Squircular mapping

$$x = \frac{r\ \ sgn(\cos t)}{\sqrt{2}\,|\sin t|}\ \sqrt{1 - |\cos 2t|}$$

$$y = \frac{r\ \ sgn(\sin t)}{\sqrt{2}\,|\cos t|}\ \sqrt{1 - |\cos 2t|}$$

## 9 Arc Length of the FG-squircle

In this section, we shall provide three ways for calculating the arc length of the FG-squircle. Unfortunately, the formulas we provide in this section are not closed-form analytical expressions. All of the arc length equations here involve integrals that have to be calculated numerically to get an approximation of arc length.

### 9.1 Rectangular Cartesian coordinates y=f(x)

Given a curve with equation y=f(x), the formula for arc length in Cartesian coordinates is

$$l\ =\ \int_{x_1}^{x_2} \sqrt{1 + [f'(x)]^2}\, dx$$

In order to use this formula for the FG-squircle, we need to convert the squircle equation to the form y=f(x). To do this, we start with the standard equation for the FG-squircle: $x^2 + y^2 - \frac{s^2}{r^2}\ x^2y^2 = r^2$
Isolate y to one side of the equation, we get

$$y^2\ =\ \frac{r^2 - x^2}{1 - \frac{s^2}{r^2}x^2}$$

$$\Rightarrow \quad f(x) = y = \sqrt{\frac{r^2 - x^2}{1 - \frac{s^2}{r^2}x^2}}$$

Differentiating with respect to x, we get

$$f'(x) = \frac{s^2 x\sqrt{r^2 - x^2}}{r^2\left(1 - \frac{s^2}{r^2}x^2\right)^{\frac{3}{2}}} - \frac{x}{\sqrt{r^2 - x^2}\sqrt{1 - \frac{s^2}{r^2}x^2}}$$

which simplifies to

$$f'(x) = \frac{r^3 x\,(s^2 - 1)}{\sqrt{r^2 - x^2}(r^2 - s^2x^2)^{\frac{3}{2}}}$$

Hence the incomplete arc length of the FG-squircle is

$$l = \int_{x_1}^{x_2} \sqrt{1 + \frac{r^6 x^2\,(s^2 - 1)^2}{(r^2 - x^2)(r^2 - s^2x^2)^3}}\, dx$$

Unfortunately, we were unable to simplify this integral further. We even tried using computer algebra software with symbolic math capabilities, but could not find a closed-form equation for the arc length of the FG-squircle. However, this integral can be calculated numerically using standard techniques in numerical analysis.

## 9.2 Parametric formula for Arc Length

Given a curve with parametric equations for coordinates x(t) and y(t), the formula for arc length is

$$l = \int_{t_1}^{t_2} \sqrt{\left(\frac{dx}{dt}\right)^2 + \left(\frac{dy}{dt}\right)^2}\, dt$$

One can then use any of the three parametric equations given in section 8 to plug into this to get an arc length formula. This integral formula can then be numerically approximated.

## 9.3 Polar coordinates $\rho = f(\theta)$

Another approach to computing the arc length of the FG-squircle is using polar coordinates. The equation for arc length in polar coordinates is:

$$l = \int_{\theta_1}^{\theta_2} \sqrt{\rho^2 + \left(\frac{d\rho}{d\theta}\right)^2}\, d\theta$$

One can then use the polar form of the FG-squircle given in section 7 and plug into this equation to get an equation for arc length which can then be numerically evaluated. For the complete arc length, simply use the limits of 0 to 2π

$$l_{complete} = \int_{0}^{2\pi} \sqrt{\rho^2 + \left(\frac{d\rho}{d\theta}\right)^2}\, d\theta$$

## 10 Three Dimensional Counterpart of the FG-squircle

### 10.1 Algebraic Surface

In this section, we present a three dimensional counterpart to the FG-squircle. The algebraic equation for this 3D shape is

$$x^2 + y^2 + z^2 - \frac{s^2}{r^2}\,x^2y^2 - \frac{s^2}{r^2}\,y^2z^2 - \frac{s^2}{r^2}\,x^2z^2 + \frac{s^4}{r^4}x^2y^2z^2 = r^2$$

Just like the FG-squircle, this shape has two parameters: squareness *s* and radius *r*. When s = 0, the shape is a sphere with radius *r*. When s =1, the shape is a cube with edge length *2r*. In between, it is a three dimensional shape that resembles the sphere and the cube. The shape is shown in Figure 8 at varying values of squareness.

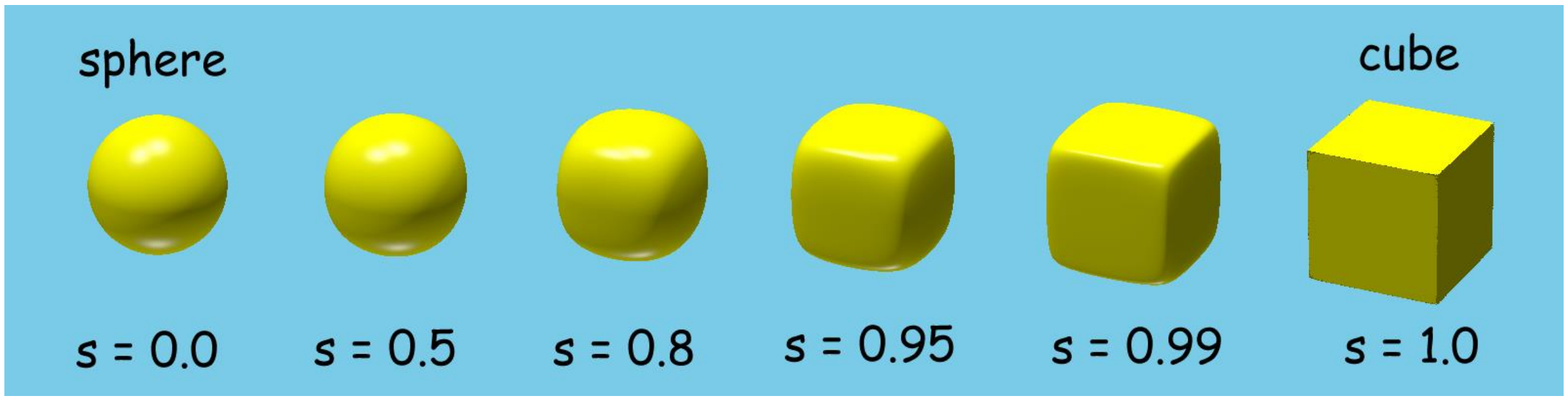

Figure 8: The 3D squircle

We propose naming this shape as the *sphube*, which is a portmanteau of sphere and cube. This sort of naming would be consistent with the naming of the squircle. This algebraic shape is a sextic surface [Weisstein 2016a] because of the $\frac{s^4}{r^4}x^2y^2z^2$ term.

The derivation of this 3D shape follows from the same logic as the 2D derivation of the FG-squircle previously given in Section 2. We start with an expression for the cube:

$$\sqrt{1 - \frac{x^2}{r^2}}\sqrt{1 - \frac{y^2}{r^2}}\sqrt{1 - \frac{z^2}{r^2}} = 0$$

Fernandez-Guasti actually came up with this equation for the cube at the end of his 1992 paper [Fernandez-Guasti 1992], but he did not carry out his analysis much further.

If we proceed by multiplying out each of the binomial terms inside the square root, we get

$$\sqrt{(1 - \frac{x^2}{r^2})(1 - \frac{y^2}{r^2})(1 - \frac{z^2}{r^2})} = 0$$

$$\Rightarrow \quad \sqrt{1 - \frac{x^2}{r^2} - \frac{y^2}{r^2} - \frac{z^2}{r^2} + \frac{x^2y^2}{r^4} + \frac{x^2z^2}{r^4} + \frac{y^2z^2}{r^4} - \frac{x^2y^2z^2}{r^6}} = 0$$

$$\Rightarrow \quad \frac{x^2}{r^2} + \frac{y^2}{r^2} + \frac{z^2}{r^2} - \frac{x^2y^2}{r^4} - \frac{x^2z^2}{r^4} - \frac{y^2z^2}{r^4} + \frac{x^2y^2z^2}{r^6} = 1$$

This is a sextic polynomial equation for the cube.

Moreover, if we introduce a squareness parametric analogous to the logic behind the FG-squircle in Section 2, we get this equation for an intermediate 3D shape between the sphere with the cube.

$$\frac{x^2}{r^2}+\frac{y^2}{r^2}+\frac{z^2}{r^2}-\frac{s^2x^2y^2}{r^4}-\frac{s^2x^2z^2}{r^4}-\frac{s^2y^2z^2}{r^4}+\frac{s^4x^2y^2z^2}{r^6}=1$$

Multiplying both sides of the equation by $r^2$ gives us the equation for the sphube.

Manuel Fernandez-Guasti already had the necessary mathematical infrastructure to carry out this 3D derivation back in 1992, but he did not. Back then, 3D computer graphics was still a nascent technology. It was not easy to visualize three-dimensional algebraic surfaces in the early '90s. Consequently, Fernandez-Guasti ended his 1992 paper with an equation for the cube and mentioned that it can be extended to higher dimensions. He hinted that the equation for the hypercube is

$$\sqrt{1-\frac{x^2}{r^2}}\sqrt{1-\frac{y^2}{r^2}}\sqrt{1-\frac{z^2}{r^2}}\sqrt{1-\frac{\omega^2}{r^2}}=0$$

Fast forward to 24 years later, computer graphics is now a mature technology and easily accessible to the masses. Moore's law has helped commoditize 3D computer graphics through cheap computing hardware. It is now quite easy to visualize three-dimensional algebraic surfaces. Consequently, we follow on Fernandez-Guasti's footsteps and introduce the three-dimensional counterpart of his squircle.

There are several techniques for visualizing algebraic surfaces using computer graphics. One way is to use ray tracing [Hanrahan 1983]. Another way is to use isosurface extraction [Bourke 1994] which converts volumetric data into a polygonal mesh. The most popular method for isosurface extraction is the Marching Cubes algorithm [Lorenson 1987].

## 10.2 Lobes of the Sphube

Just like in the two-dimensional FG-squircle shown in Figure 2, the sphube has disjoint portions outside of the central core. For our purposes, these extraneous surfaces are unwanted and need to be pruned away in order to effectively study and visualize the central core. There are 12 extraneous lobes corresponding to each of the edges of the cube. Unlike the bounded central core, these lobes have infinite span.

Figure 9 shows the sphube along with its lobes at several different squareness values. To highlight the distinction between the central core surface and its lobes, we have colored the central core as yellow and the surrounding lobes as green. Note that only 8 of the lobes are visible because 4 of front lobes have been removed in the diagrams. This was done because the 4 lobes in front obscure the central core in the 3D view.

It is possible to cull out all the unwanted lobes to get the central squircular core shown in Figure 8. This can be done by simply restricting the surface points (x,y,z) to these inequalities:

$$-r\le x\le r$$

$$-r\le y\le r$$

$$-r\le z\le r$$

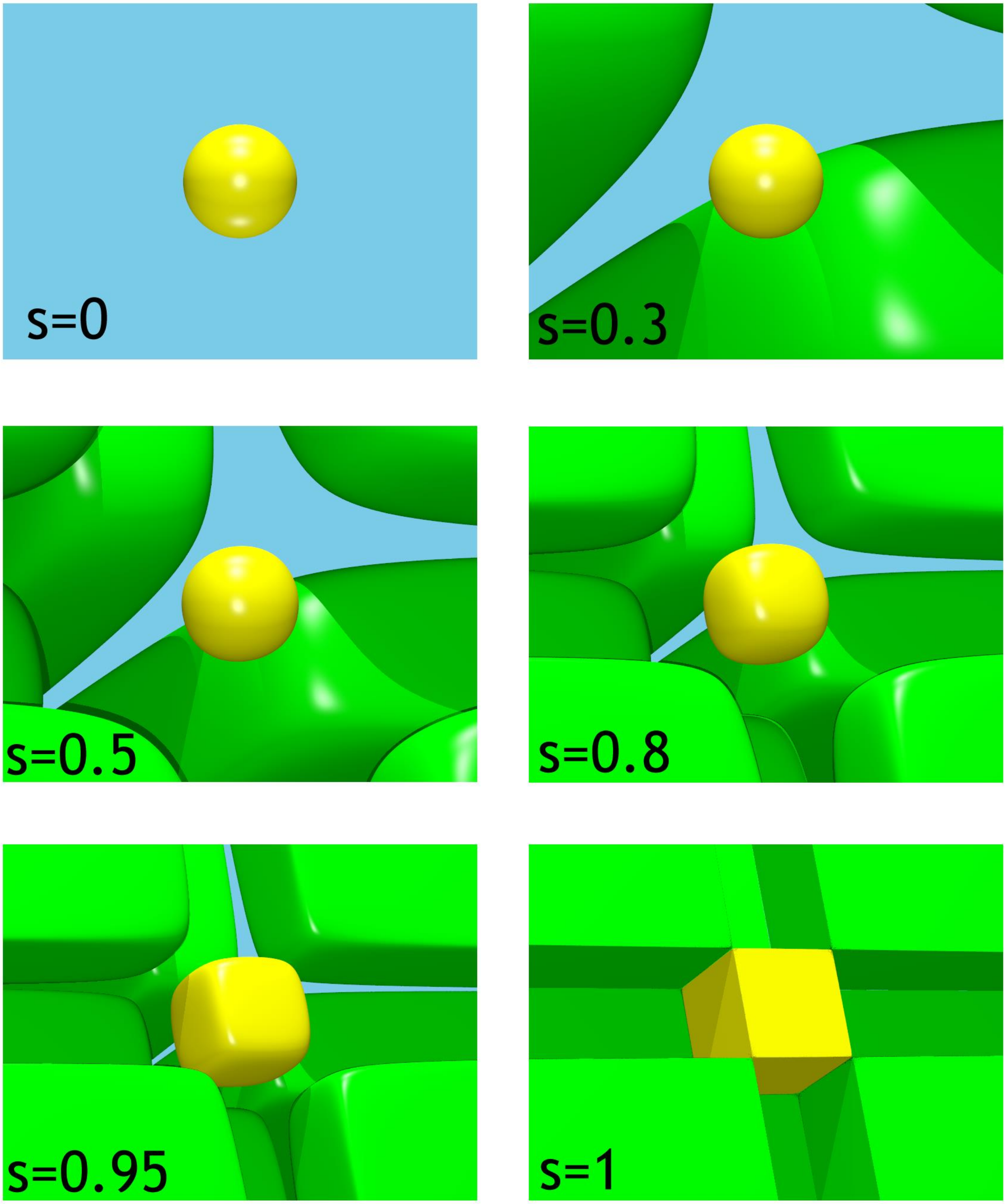

Figure 9: The sphube with lobes at different values of the squareness parameter
An animated video of the sphube is available in: https://youtu.be/f1nN9DfhgbU

# 11 Higher Dimensional Counterparts

## 11.1 Four Dimensional FG-Squircle

In this section, we will discuss higher dimensional counterparts of the FG-squircle. First, we would start with the equation for the four dimensional counterpart of the FG-squircle. We shall use the variables $x$, $y$, $z$, and $\omega$ (omega) to represent the hypersurface in the four dimensional space $\mathbb{R}^4$.

The four dimensional version of the sphere is known as the hypersphere (also 3-sphere). Similarly, the four dimensional version of the cube is known as the hypercube (also tesseract). We believe that there is an algebraic hypersurface that is an intermediate shape between these two four-dimensional shapes. In fact, we can extend the derivation given in Section 2 to four dimensions, just as we did previously for three dimensions. We propose naming this hypersurface as the *hypersphube*.

By carrying out the same logical expansion as in the previous section, we can deduce that the four dimensional counterpart of the FG-squircle has this algebraic equation:

$$x^2 + y^2 + z^2 + \omega^2 - \frac{s^2}{r^2}\,x^2y^2 - \frac{s^2}{r^2}\,y^2z^2 - \frac{s^2}{r^2}\,x^2z^2 - \frac{s^2}{r^2}\,x^2\omega^2 - \frac{s^2}{r^2}\,y^2\omega^2 - \frac{s^2}{r^2}\,z^2\omega^2$$

$$+\frac{s^4}{r^4}x^2y^2z^2 + \frac{s^4}{r^4}x^2y^2\omega^2 + \frac{s^4}{r^4}x^2z^2\omega^2 + \frac{s^4}{r^4}y^2z^2\omega^2 - \frac{s^6}{r^6}x^2y^2z^2\omega^2 = r^2$$

This algebraic equation has order-8 (octic) and has 16 terms. The octic term $-\frac{s^6}{r^6}x^2y^2z^2\omega^2$ has the highest effective exponent in the equation. In general, an n-dimensional counterpart of the FG-squircle would have $2^n$ terms and have a term with the highest exponent of order $2n$. We will discuss the n-dimensional counterpart of the FG-squircle in more detail in the next subsection.

There have been many proposals on how to visualize four-dimensional shapes [Banchoff 1990]. Techniques include using projections, 3D cross sections, and colored animations. We shall relegate the task of visualizing and studying the hypersphube as future work.

## 11.2 Higher dimensional FG-Squircle

Using the same reasoning as the three dimensional counterpart of the FG-Squircle, it is possible to get an n-dimensional analogue. We assert without proof that the equation for the n-dimensional FG-Squircle is

$$0 = \sum_{j=1}^{2^n} (-1)^{\psi(j)+1}\ \frac{s^{2\psi(j)-2}}{r^{2\psi(j)}} \prod_{i=1}^{n} x_i^{2\alpha(i,j)}$$

This equation is quite complicated and will require some definitions for the terms. In particular we need to define the variables $x_i$ and the functions $\psi(j)$ and $\alpha(i,j)$ **.** The dummy variables $i$ and $j$ are just indices used in the summation and product terms of the formula. Here are the definitions:

- $x_i$ is the *ith* variable representing an axis of *n*-dimensional space. For example, in 3-D, we have $x_1$ as a stand-in for $x$; $x_2$ as a stand-in for $y$; and $x_3$ as a stand-in for $z$. In higher dimensions, we run out of letters for variables in each axis, so we have to resort to using subscripts instead.
- $\psi(j)$ is the popcount function. The name is short for population count. It is the sum total of the number of 1-bits in the binary representation of the number $j-1$.
- $\alpha(i,j)$ is a bit extraction function. It extracts the *ith* bit in the binary representation of the number $j-1$. Note that $\alpha(i,j)$ is a binary function; i.e. it can only produce a value of $0$ or $1$.

The popcount function [Warren 2013] is quite famous and useful in the field of Computer Science. In fact, most modern microprocessors support it in hardware. All Intel processors with SSE4.2 instruction set extension have a built-in instruction for calculating popcount.

# 12 Mapping a Cube to a Sphere

In this section, we shall present a way to map the solid cube to a solid sphere. This would be a volumetric mapping which takes a point inside of the cube and map it to a unique point inside the sphere. In the next section, we shall provide the inverse of this mapping. Before we start, we need to need to define the canonical space for the mapping. We shall define our sphere as the set $S = \{(u, v, w) |\ u^2 + v^2 + w^2 \leq 1\}$. In other words, our sphere will be the unit sphere centered at the origin. Similarly, we shall define our cube as the set $C = [-1,1]$ x $[-1,1$x $[-1,1]$. This cube is also centered at the origin with a side length of 2. We shall denote *(u,v,w)* as points inside the sphere, and *(x,y,z)* as points inside the cube.

In analogy to the 2D case, we shall use the sphube for mapping from the cube to a sphere. The idea is a direct 3D extension of the derivation of the FG-squircular mapping in [Fong 2014]. We can think of the unit solid sphere as continuum of concentric hollow spheres growing in radius. Likewise, we can think of the solid cube as a continuum of concentric sphubes growing in radius.

Also in analogy to the 2D case [Fong 2014], we can introduce the concept of the *shrunken sphube*. This is sphube with the given constraint of s = r in its equation. Basically, we have a parameter s that controls both the squareness and the size of the sphube. Visually, we have a continuum of growing shapes that start as a point at s=0 and end as a cube at s=1. At s = ½, we have a half-sized sphube that is neither spherical nor cubic in shape, but in between the two. The equation for the shrunken sphube is

$$x^2 + y^2 + z^2 - x^2y^2 - y^2z^2 - x^2z^2 + x^2y^2z^2 = s^2$$

In yet another analogy to the 2D case [Fong 2014], we can introduce the concept of the *squircularity condition*. This is basically an equation assigning each hollow sphere inside the solid sphere to a shrunken sphube inside the solid cube. The equation for this constraint is

$$u^2 + v^2 + w^2 = x^2 + y^2 + z^2 - x^2y^2 - y^2z^2 - x^2z^2 + x^2y^2z^2$$

Using the squircularity condition and a constraint that our mapping be radial in nature, we can get the following equations for the mapping:

$$u = x\frac{\sqrt{x^2 + y^2 + z^2 - x^2y^2 - y^2z^2 - x^2z^2 + x^2y^2z^2}}{\sqrt{x^2 + y^2 + z^2}}$$

$$v = y\frac{\sqrt{x^2 + y^2 + z^2 - x^2y^2 - y^2z^2 - x^2z^2 + x^2y^2z^2}}{\sqrt{x^2 + y^2 + z^2}}$$

$$w = z\frac{\sqrt{x^2 + y^2 + z^2 - x^2y^2 - y^2z^2 - x^2z^2 + x^2y^2z^2}}{\sqrt{x^2 + y^2 + z^2}}$$

In order to avoid a division by zero at the origin where *(x,y,z) = (0,0,0),* we shall impose that *u=0, v=0,* and *w=0* when the input point to the mapping is the origin. Note that when exactly one of the coordinates is zero, the equations reduce to a 2D case which is equivalent to the FG-squircular mapping.

# 13 Mapping a Sphere to a Cube

We shall now invert the volumetric mapping from the previous section. Before, we proceed we need to write down the equations governing the radial constraint on the mapping. As a direct analogy to the 2D case [Fong 2014] and using symmetry, the following equations must hold in order for the mapping to be radial.

$$y = \frac{v}{u}x \qquad z = \frac{w}{u}x \qquad y = \frac{v}{w}z$$

Using the squircularity condition introduced in the previous section and the radial equations from above, we can derive a polynomial equation in *x* that we need for computing the inverse.

$$u^2 + v^2 + w^2 = x^2 + y^2 + z^2 - x^2y^2 - y^2z^2 - x^2z^2 + x^2y^2z^2$$

$$= x^2 + \frac{v^2}{u^2}x^2 + \frac{w^2}{u^2}x^2 - x^2\frac{v^2}{u^2}x^2 - \frac{v^2}{u^2}x^2\frac{w^2}{u^2}x^2 - x^2\frac{w^2}{u^2}x^2 + x^2\frac{v^2}{u^2}x^2\frac{w^2}{u^2}x^2$$

$$= x^2 + \frac{v^2}{u^2}x^2 + \frac{w^2}{u^2}x^2 - \frac{v^2}{u^2}x^4 - \frac{v^2}{u^2}\frac{w^2}{u^2}x^4 - \frac{w^2}{u^2}x^4 + \frac{v^2}{u^2}\frac{w^2}{u^2}x^6$$

$$= \frac{v^2w^2}{u^4}x^6 - \left(\frac{v^2}{u^2} + \frac{w^2}{u^2} + \frac{v^2w^2}{u^4}\right)x^4 + \left(1 + \frac{v^2}{u^2} + \frac{w^2}{u^2}\right)x^2$$

We can then move all the terms to one side of the equation to get this sextic polynomial equation in x.

$$\frac{v^2w^2}{u^4}x^6 - \left(\frac{v^2}{u^2} + \frac{w^2}{u^2} + \frac{v^2w^2}{u^4}\right)x^4 + \left(1 + \frac{v^2}{u^2} + \frac{w^2}{u^2}\right)x^2 - (u^2 + v^2 + w^2) = 0$$

Fortunately, this sextic equation has zero coefficients for all the odd powers. In other words, we do not have $x^5$, $x^3$, and linear x terms in the polynomial equation. So we can solve for $x^2$ using the standard cubic formula, and then take the square root to solve for x. Wikipedia provides a good explanation and reference for the cubic formula. We will use the following coefficients as input to the cubic formula.

$$a = \frac{v^2w^2}{u^4}$$

$$b = -\left(\frac{v^2}{u^2} + \frac{w^2}{u^2} + \frac{v^2w^2}{u^4}\right)$$

$$c = 1 + \frac{v^2}{u^2} + \frac{w^2}{u^2}$$

$$d = -(u^2 + v^2 + w^2)$$

The standard cubic formula produces unwieldy equations in terms of the 4 coefficients given above. We shall therefore introduce intermediate variables to make the equations simpler. First, we need the discriminant for the cubic formula, which we will denote as the intermediate variable $D_o$

$$D_o = 18abcd - 4b^3d + b^2c^2 - 4ac^3 - 27a^2d^2$$

Next, we introduce another intermediate variable $C_o$. It is defined as a cube root expression

$$C_o = \sqrt[3]{\frac{2b^3 - 9abc + 27a^2d + \sqrt{-27a^2D_o}}{2}}$$

Armed with the intermediate variables $D_o$ and $C_o$, we use the standard cubic formula to get an equation for $x^2$.

$$x^2 = -\frac{1}{3a}(b + C_o + \frac{b^2 - 3ac}{C_o})$$

We can then take the square root of this equation to get x. Furthermore, in order to account for the correct sign of x, we shall use the signum function.

$$x = sgn(u)\sqrt{-\frac{1}{3a}(b + C_o + \frac{b^2 - 3ac}{C_o})}$$

Note that in this formula, the intermediate variables a, b ,c, and $C_o$ can all be calculated from u,v, and w. Albeit, expanding the equation for x in terms of u, v, and w will produce a very unwieldy and long algebraic equation. So, we will leave the expression for x in its current form.

The equations for y and z are pretty simple once we have calculated x. These are

$$y = \frac{v}{u}x \qquad\qquad z = \frac{w}{u}x$$

Note that the equations provided above assume that none of u, v, or w are zero. We have not singled out these special cases above. Special care must use to handle such cases. The rules are:

1) Case 1: when u = v = w = 0
   In this case, we are at the origin. Simply set (x,y,z) = (0,0,0)

2) Case 2: exactly 2 of u, v, and w are zero.
   Assign x=u, y=v, and w=z

3) Case 3: exactly 1 of u, v, or w is zero
   This case is tricky, but it reduces the 3D mapping to a 2D mapping. We can just use an analogous disc-to-square equation in the FG-squircular mapping shown in Figure 5.

Also note that when the discriminant $D_o > 0$, the expression for $C_o$ will involve the cube root of a complex number. In this situation, De Moivre's formula can be used for calculating the complex-valued $C_o$ intermediate variable. If the software implementation can handle complex numbers, the given equation above for x is valid and will produce the correct real number results even though there are intermediate complex numbers. If one wishes to avoid the use of complex numbers altogether in the formulas, extra work must be done.

To avoid using complex numbers in the inverse equations, we need to introduce more intermediate variables and steps in calculating the inverse. These extra equations only hold for $D_o>0$.

$$r_o = \frac{2b^3 - 9abc + 27a^2d}{2}$$

$$q_o = \tfrac{1}{2}\sqrt{27a^2D_o}$$

Both the intermediate variables $r_o$ and $q_o$ are real numbers that are components of the complex number radicand inside the cube root expression for $C_o$.

Using De Moivre's formula on $C_o$, the expression for $x^2$ simplifies to

$$x^2 = -\frac{1}{3a}\left( b + 2\,(r_o^2 + q_o^2)^{\frac{1}{6}} \cos\frac{\tan^{-1}\frac{q_o}{r_o}}{3} \right)$$

This expression only deals with intermediate variables that are real numbers. Furthermore, taking the square root of this expression gives the desired value of x.

In summary, here are the equations for solving the inverse:

if $D_o \leq 0$

$$\boldsymbol{x = sgn(u)\sqrt{-\frac{1}{3a}\left(b + C_o + \frac{b^2 - 3ac}{C_o}\right)}}$$

else

$$\boldsymbol{x = sgn(u)\sqrt{-\frac{1}{3a}\left( b + 2\,(r_o^2 + q_o^2)^{\frac{1}{6}} \cos\frac{\tan^{-1}\frac{q_o}{r_o}}{3} \right)}}$$

The equations for y and z remain the same.

$$\boldsymbol{y = \frac{v}{u}x}$$

$$\boldsymbol{z = \frac{w}{u}x}$$

All these equations are derived using the standard cubic formula and De Moivre's formula which are described in Wikipedia.

Since the equations for the mappings are not simple, we have included a reference software implementation of the mappings in the paper appendix.

## 14 Three Dimensional Squelched Mapping

In this section, we shall come up with a three dimensional analogue of the squelched grid mapping [Fong 2017]. This would be a volumetric mapping that maps a point inside of the sphere to a point inside the cube.

### 14.1 Sphere to Cube

Just like the two dimensional case of squelched grid mapping starts with 2D ellipse constraint equations, we have 3D ellipsoid constraint equations. These are

$$1 = \frac{u^2}{x^2} + v^2 + w^2$$
$$1 = u^2 + \frac{v^2}{y^2} + w^2$$
$$1 = u^2 + v^2 + \frac{w^2}{z^2}$$

We can then do some simple algebraic manipulations to isolate x, y, and z from each of these three equations, whereupon we get these sphere to cube equations

$$\boldsymbol{x = \frac{u}{\sqrt{1 - v^2 - w^2}}}$$

$$\boldsymbol{y = \frac{v}{\sqrt{1 - u^2 - w^2}}}$$

$$\boldsymbol{z = \frac{w}{\sqrt{1 - u^2 - v^2}}}$$

### 14.2 Cube to Sphere

Here, we invert the sphere-to-cube equations above to get cube-to-sphere equations. First, we start with the given equations and do some algebraic manipulation

$$x = \frac{u}{\sqrt{1 - v^2 - w^2}} \quad \Rightarrow \quad u^2 = x^2(1 - v^2 - w^2)$$

$$y = \frac{v}{\sqrt{1 - u^2 - w^2}} \quad \Rightarrow \quad v^2 = y^2(1 - u^2 - w^2)$$

$$z = \frac{w}{\sqrt{1 - u^2 - v^2}} \quad \Rightarrow \quad w^2 = z^2(1 - u^2 - v^2)$$

Substitute the value for *w* into the first 3D ellipsoid constraint equation to get

$$1 = \frac{u^2}{x^2} + v^2 + w^2 = \frac{u^2}{x^2} + v^2 + z^2(1 - u^2 - v^2) = \frac{u^2}{x^2} + v^2 + z^2 - u^2z^2 - v^2z^2$$

$$\Rightarrow \quad 1 - v^2 - z^2 + v^2z^2 = \frac{u^2}{x^2} - u^2z^2 \quad \Rightarrow \quad u^2 = x^2\,\frac{1 - v^2 - z^2 + v^2z^2}{1 - x^2z^2}$$

Similarly, we can substitute the value for w into the second 3D ellipsoid constraint equation to get

$$1 = u^2 + \frac{v^2}{y^2} + w^2 = u^2 + \frac{v^2}{y^2} + z^2(1 - u^2 - v^2) = u^2 + \frac{v^2}{y^2} + z^2 - u^2z^2 - v^2z^2$$

$$\Rightarrow \quad 1 - u^2 - z^2 + u^2z^2 = \frac{v^2}{y^2} - v^2z^2 \quad \Rightarrow \quad v^2 = y^2\,\frac{1 - u^2 - z^2 + u^2z^2}{1 - y^2z^2}$$

whereupon we can substitute this expression for $v^2$ into the expression for $u^2$ to get

$$u^2 = x^2\ \frac{1 - v^2 - z^2 + v^2 z^2}{1 - x^2 z^2} \ = \ x^2\ \frac{1 - y^2\ \dfrac{1 - u^2 - z^2 + u^2 z^2}{1 - y^2 z^2} - z^2 + z^2 y^2\ \dfrac{1 - u^2 - z^2 + u^2 z^2}{1 - y^2 z^2}}{1 - x^2 z^2}$$

$$\Rightarrow \qquad u^2 = x^2 \frac{1 - y^2 z^2 - y^2(1 - u^2 - z^2 + u^2 z^2) - z^2(1 - y^2 z^2) + y^2 z^2(1 - u^2 - z^2 + u^2 z^2)}{(1 - x^2 z^2)(1 - y^2 z^2)}$$

$$\Rightarrow \qquad u^2 = x^2 \frac{1 - y^2 - z^2 + y^2 z^2 + u^2 y^2(1 - 2z^2 + z^4)}{(1 - x^2 z^2)(1 - y^2 z^2)}$$

We can then isolate $u^2$ into one side of the equation

$$u^2 - u^2 x^2 y^2 \frac{1 - 2z^2 + z^4}{(1 - x^2 z^2)(1 - y^2 z^2)} = x^2 \frac{1 - y^2 - z^2 + y^2 z^2}{(1 - x^2 z^2)(1 - y^2 z^2)}$$

$$\Rightarrow \qquad u^2 \frac{(1 - x^2 z^2)(1 - y^2 z^2) - x^2 y^2(1 - 2z^2 + z^4)}{(1 - x^2 z^2)(1 - y^2 z^2)} = x^2 \frac{1 - y^2 - z^2 + y^2 z^2}{(1 - x^2 z^2)(1 - y^2 z^2)}$$

$$\Rightarrow \qquad u^2 \frac{1 - x^2 y^2 - y^2 z^2 - x^2 z^2 + 2x^2 y^2 z^2}{(1 - x^2 z^2)(1 - y^2 z^2)} = x^2 \frac{(1 - y^2)(1 - z^2)}{(1 - x^2 z^2)(1 - y^2 z^2)}$$

$$\Rightarrow \qquad u^2 = x^2 \frac{(1 - y^2)(1 - z^2)}{1 - x^2 y^2 - y^2 z^2 - x^2 z^2 + 2x^2 y^2 z^2}$$

This gives us an expression for *u* in terms of *x*, *y*, and *z*.

$$\boldsymbol{u = x\sqrt{\frac{(1 - y^2)(1 - z^2)}{1 - x^2 y^2 - y^2 z^2 - x^2 z^2 + 2x^2 y^2 z^2}}}$$

Similarly, one can get these expressions for *v* and *w*.

$$\boldsymbol{v = y\sqrt{\frac{(1 - x^2)(1 - z^2)}{1 - x^2 y^2 - y^2 z^2 - x^2 z^2 + 2x^2 y^2 z^2}}}$$

$$\boldsymbol{w = z\sqrt{\frac{(1 - x^2)(1 - y^2)}{1 - x^2 y^2 - y^2 z^2 - x^2 z^2 + 2x^2 y^2 z^2}}}$$

### 14.3 Parametric Equations for the Sphube

Using the squelched mapping equations for sphere-to-cube, we can come up with parametric equations for the sphube surface discussed in section 10. We shall proceed in analogy to the 2D case in Section 8 where we made the substitution $u = cos(t)$ and $v = sin(t)$ to come up with parametric equations for the 2D squircle. This time we will do a substitution using the canonical parametric equations for the sphere. The equations for these are

$$u = r\cos\theta\cos\phi$$
$$v = r\cos\theta\sin\phi$$
$$w = r\sin\theta$$

We can then plug-in these values for *u,v*, and *w* into the sphere-to-cube equations in Section 14.1 and simplify. This will give us parametric equations for the sphube

$$\boldsymbol{x = \frac{r\,cos\,\theta\,cos\,\phi}{\sqrt{1 - s_p\,cos^2\,\theta\,sin^2\,\phi - s_p\,sin^2\,\theta}}}$$

$$\boldsymbol{y = \frac{r\,cos\,\theta\,sin\,\phi}{\sqrt{1 - s_p\,cos^2\,\theta\,cos^2\,\phi - s_p\,sin^2\,\theta}}}$$

$$\boldsymbol{z = \frac{r\,sin\,\theta}{\sqrt{1 - s_p\,cos^2\,\theta}}}$$

Here $s_p$ is a squareness value for the parametric equation, but it is not identical to $s$. Nevertheless, when $s_p$ is zero, the surface is a sphere. Likewise, when $s_p$ approaches one, the surface will appear more and more like a cube. One can easily verify these parametric equations using some 3D plotting software. For example, one can use the ParametricPlot3D command in Mathematica to visualize this surface at different values of $s_p$.

# 15 Squircular Quadrics

Quadrics are 3D surfaces with a quadratic polynomial equation. Most quadric surfaces have a circle or an ellipse as a cross section of the surface. In this section of the paper, we shall present a generalization of quadrics by using the FG-squircle in lieu of the circle. We shall call these surfaces as *squircular quadrics* (or *squadrics* for short)

Before we start with the 3D squircular shapes, we will briefly review a generalization of the FG-squircle known as the *rectellipse* [Weisstein 2016c]. This 2D shape is a stretched version of the FG-squircle where the length of the major axis is different from the length of the minor axis. The equation for the rectellipse is

$$\frac{x^2}{a^2}+\frac{y^2}{b^2}-\frac{s^2x^2y^2}{a^2b^2}=1$$

When s=0, the rectellipse is an ellipse. When *s=1*, the rectellipse is a rectangle. In between, it is an intermediate shape between the two.

Note that in this section, we will reuse the variables *a*, *b*, and *c* in a completely different context as previously used in Section 13. In the context of quadrics, it is common convention to use these variables to denote semi-major and semi-minor lengths. Hopefully, since this section is self-contained and independent of the cubic equation in Section 13, there will be no confusion or ambiguity.

We will also use a convention where the main axes of our shapes lie on xy-plane. This means that z direction is normal to the squircular cross sections. By symmetry, we can interchange the x, y, z variables in the equations to get a rotated version of the shapes.

## 15.1 Squircular Ellipsoid

The equation for the squircular ellipsoid is very similar to the equation for the sphube, except that we introduce eccentricity through asymmetrical variables *a*, *b*, and *c* denoting semi-major and semi-minor lengths. Just like the sphube, the equation for this shape is a sextic polynomial. Also, there are additional bounds appearing as inequalities for x, y, and z in order to prune away unwanted lobes.

$$\frac{x^2}{a^2}+\frac{y^2}{b^2}+\frac{z^2}{c^2}-\frac{s^2x^2y^2}{a^2b^2}-\frac{s^2x^2z^2}{a^2c^2}-\frac{s^2y^2z^2}{b^2c^2}+\frac{s^4x^2y^2z^2}{a^2b^2c^2}=1 \qquad with \qquad \begin{array}{c} -a\le x\le a \\ -b\le y\le b \\ -c\le z\le c \end{array}$$

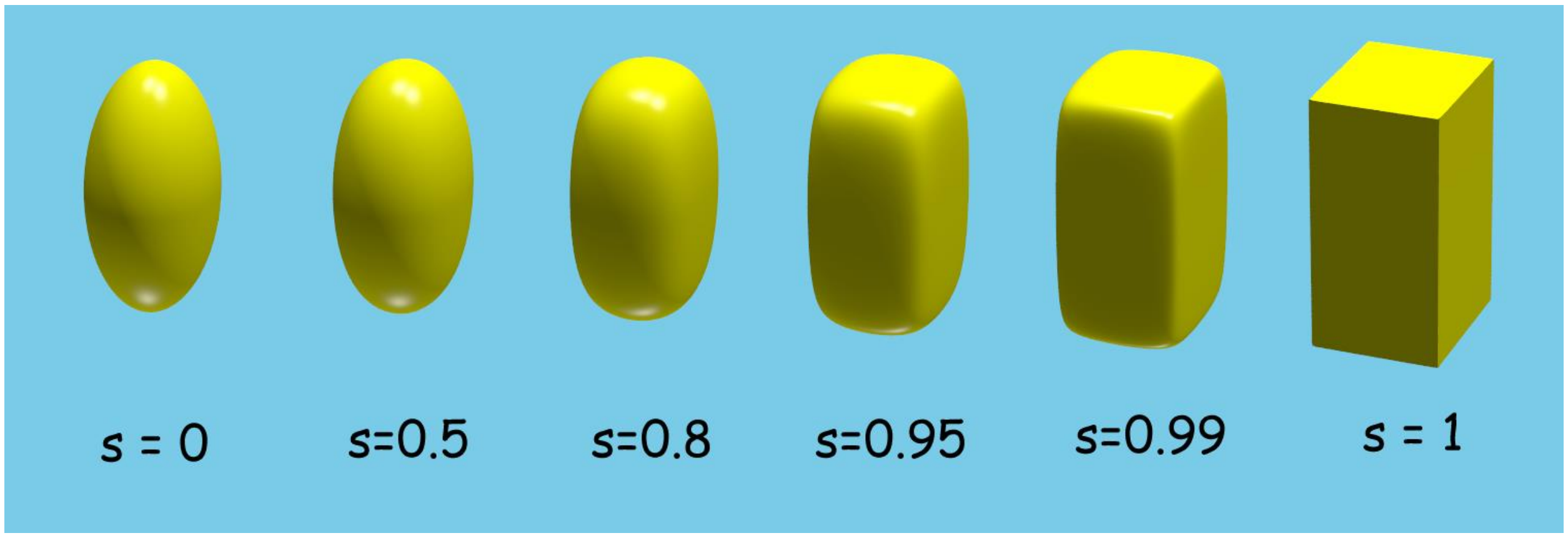

Figure 10: The squircular ellipsoid at varying squareness
with a=1.2, b=1.0, and c=2.0

### 15.2 Squircular Cylinder (*sqylinder*)

The sqylinder is simply an extrusion of the FG-squircle to three dimensions. Its equation is quite simple and does not vary with the value of z. The equation for this shape is a quartic polynomial. There are additional bounds appearing as inequalities for x, y, and z in order to prune away unwanted lobes.

$$\frac{x^2}{a^2} + \frac{y^2}{b^2} - \frac{s^2x^2y^2}{a^2b^2} = 1 \qquad with \qquad \begin{array}{r} -a \le x \le a \\ -b \le y \le b \\ 0 \le z \le c \end{array}$$

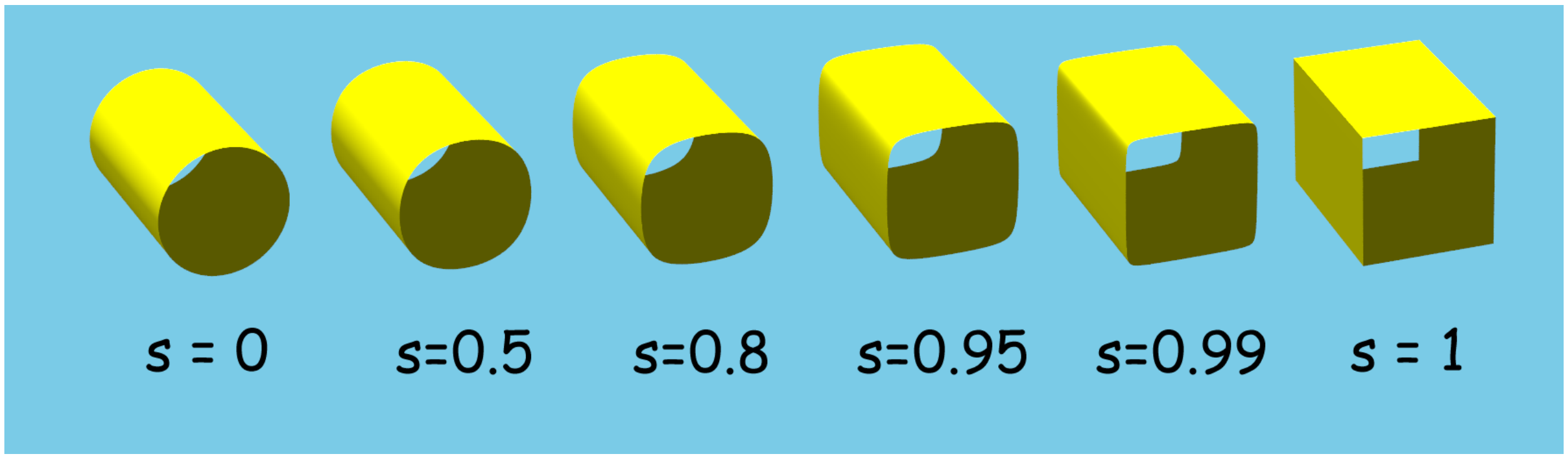

Figure 11: Hollowed squircular cylinder at varying squareness values
with a=b=1.0 and c=3.0

### 15.3 Non-Uniform Squircular Cylinder

The shape is similar to the cylinder except that it does not have a constant cross section. The cross section is modulated by the squareness parameter of the FG-squircle. At one end, it has a circular cross section. At the other end, it has a square cross section. This cylindrical shape is hollow with no end caps. The equation for this shape is a sextic polynomial. There are additional bounds appearing as inequalities for x, y, and z in order to prune away unwanted lobes.

$$\frac{x^2}{a^2} + \frac{y^2}{b^2} - \frac{z^2x^2y^2}{a^2b^2c^2} = 1 \qquad with \qquad \begin{array}{r} -a \le x \le a \\ -b \le y \le b \\ 0 \le z \le c \end{array}$$

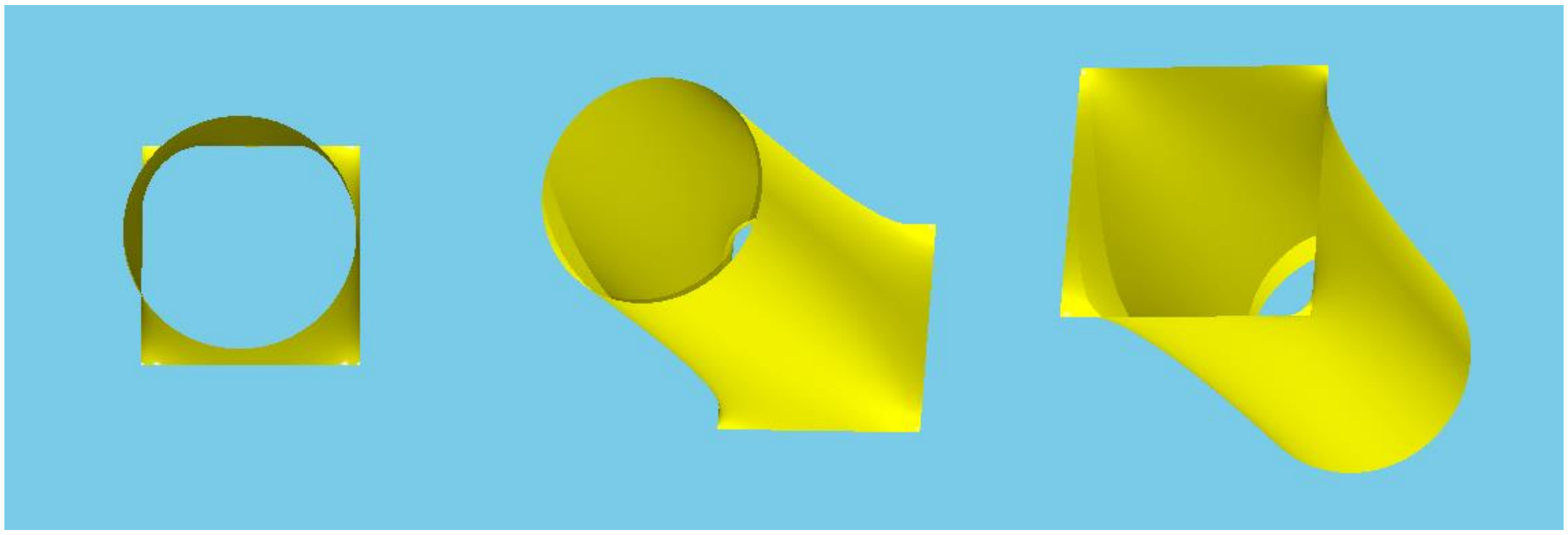

Figure 12: Top and side views of the non-uniform squircular cylinder
with a=b=1.0 and c=4.0

## 15.4 Squircular Cone

The squircular cone is an intermediate shape between the cone and the pyramid. The radius of the squircular cross section is modulated by the z value. The equation for this shape is a quartic polynomial.

$$x^2 + y^2 - \frac{s^2c^2x^2y^2}{z^2} = \frac{z^2}{c^2}$$

This equation simplifies to

$$x^2z^2 + y^2z^2 - s^2c^2x^2y^2 - \frac{z^4}{c^2} = 0$$

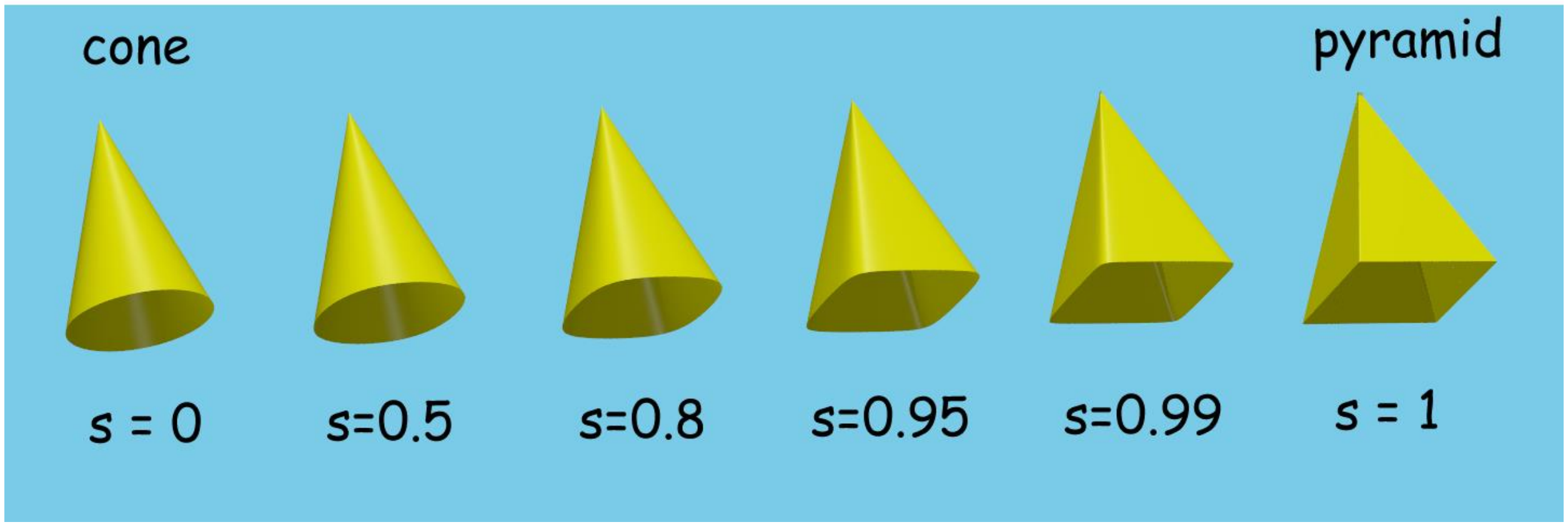

Figure 13: The squircular cone at varying squareness values with c=3.0

In addition to the implicit equation provided above, we need to subject the squircular cone to some inequalities that bound x, y, and z in order to prune away unwanted lobes in the shape.

$$0 \le z \le c \qquad -\frac{z}{c} \le x \le \frac{z}{c} \qquad -\frac{z}{c} \le y \le \frac{z}{c}$$

## 15.5 Non-Uniform Squircular Cone

The shape is similar to the cone except that it does not have a fixed-shape cross section. The cross sectional squircle varies in size and squareness with z.  At the tip on top, the cross section is circular. At the base on the bottom, the cross section is a square. The equation for this shape is a quartic polynomial. There are additional bounds appearing as inequalities for x, y, and z in order to prune away unwanted lobes.

$$\frac{x^2}{a^2} + \frac{y^2}{b^2} - \frac{x^2y^2}{a^2b^2} = \frac{z^2}{c^2} \qquad with \qquad \begin{array}{r} -a \le x \le a \\ -b \le y \le b \\ 0 \le z \le c \end{array}$$

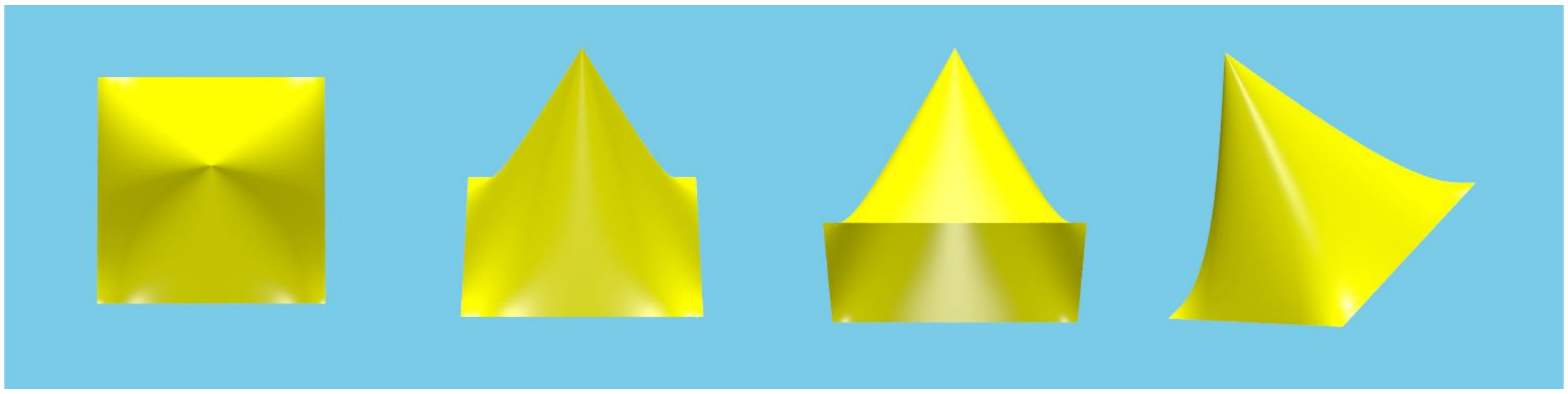

Figure 14: Different views of the non-uniform squircular cone with a=b=1.0 and c=2.0

### 15.6 Squircular Cushion

The shape resembles a pillow. It has a circular cross section on two of its main axes and a squircular cross section on remaining one. When *s=0*, this shape reverts to a sphere. The equation for this shape is a quartic polynomial

$$\frac{x^2}{a^2}+\frac{y^2}{b^2}+\frac{z^2}{c^2}-\frac{s^2x^2y^2}{a^2b^2}=1 \qquad with \qquad \begin{matrix} -a\le x\le a \\ -b\le y\le b \\ -c\le z\le c \end{matrix}$$

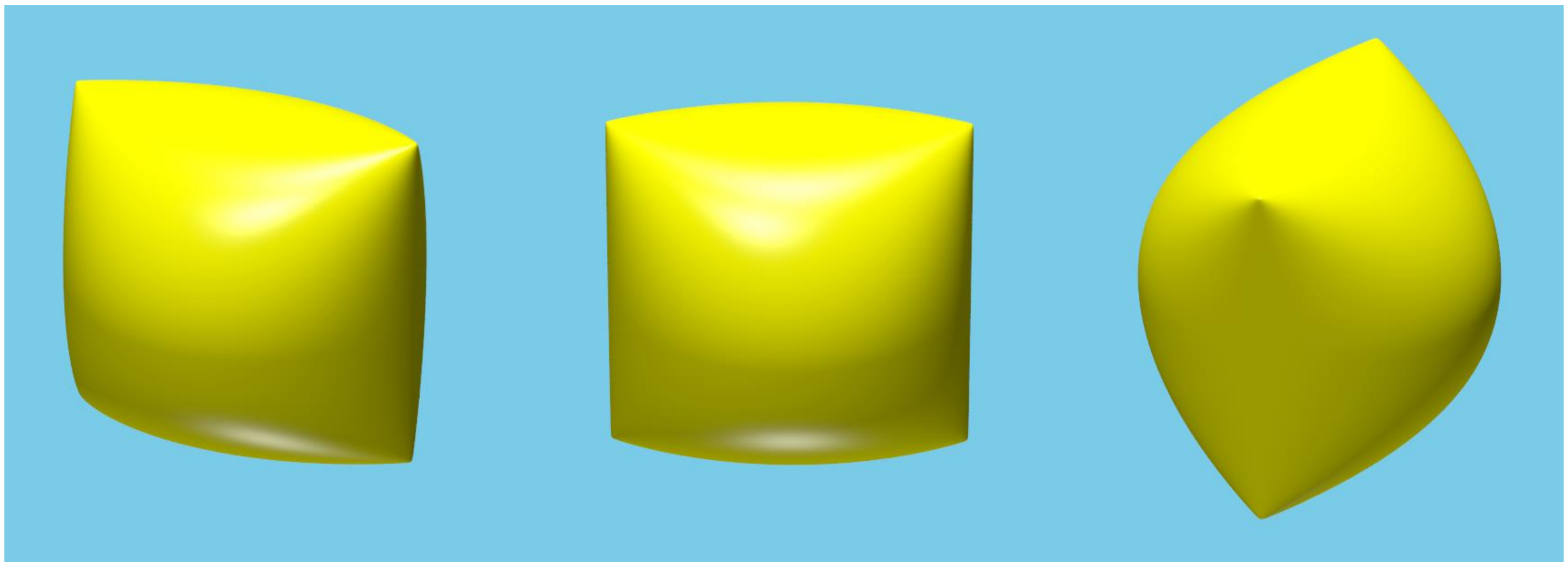

Figure 15: Different views of the squircular cushion
with a=b=c=1.0 and s=1.0

### 15.7 Sphylinder

The sphylinder is an intermediate shape between the sphere and the cylinder. Unlike the standard cylindrical surface taught in Calculus classes with no end caps and infinite length, our surface includes the two end caps of the cylinder. When *s=0*, the sphylinder reverts to a sphere with radius *r*. When *s=1*, the sphylinder becomes a cylinder with a height of *2r*. In between, we get a shape that resembles both the sphere and the cylinder. The equation for the sphylinder is a quartic polynomial.

$$x^2+y^2+z^2-\frac{s^2}{r^2}\,x^2z^2-\frac{s^2}{r^2}\,y^2z^2=r^2 \qquad with \qquad \begin{matrix} x^2+y^2\le r^2 \\ -r\le z\le r \end{matrix}$$

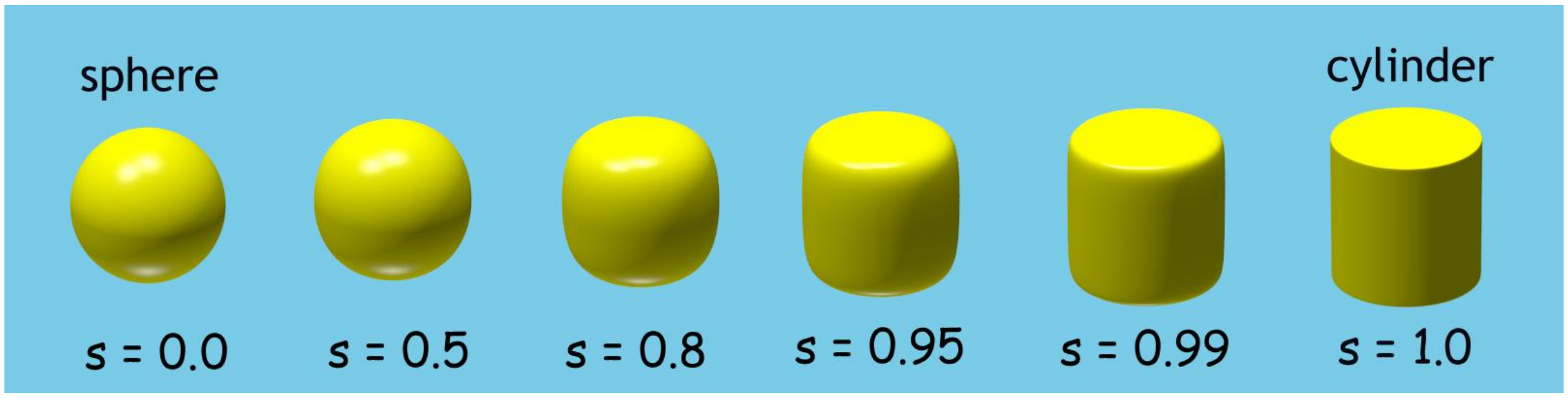

Figure 16: The sphylinder at varying squareness values
Note that the end caps are an intrinsic part of the surface.

The sphylinder is a more useful shape if we introduce eccentricity. We can introduce eccentric coefficients $a$, $b$, and $c$ to make the sphylinder have variable height and an oval-shaped cross section. However, after introducing eccentricity, the resulting shape at $s=0$ will no longer be a sphere but will be an ellipsoid instead. The equation for the eccentric sphylinder is a quartic polynomial.

$$\frac{x^2}{a^2}+\frac{y^2}{b^2}+\frac{z^2}{c^2}-\frac{s^2}{a^2c^2}\,x^2z^2-\frac{s^2}{b^2c^2}\,y^2z^2=1$$

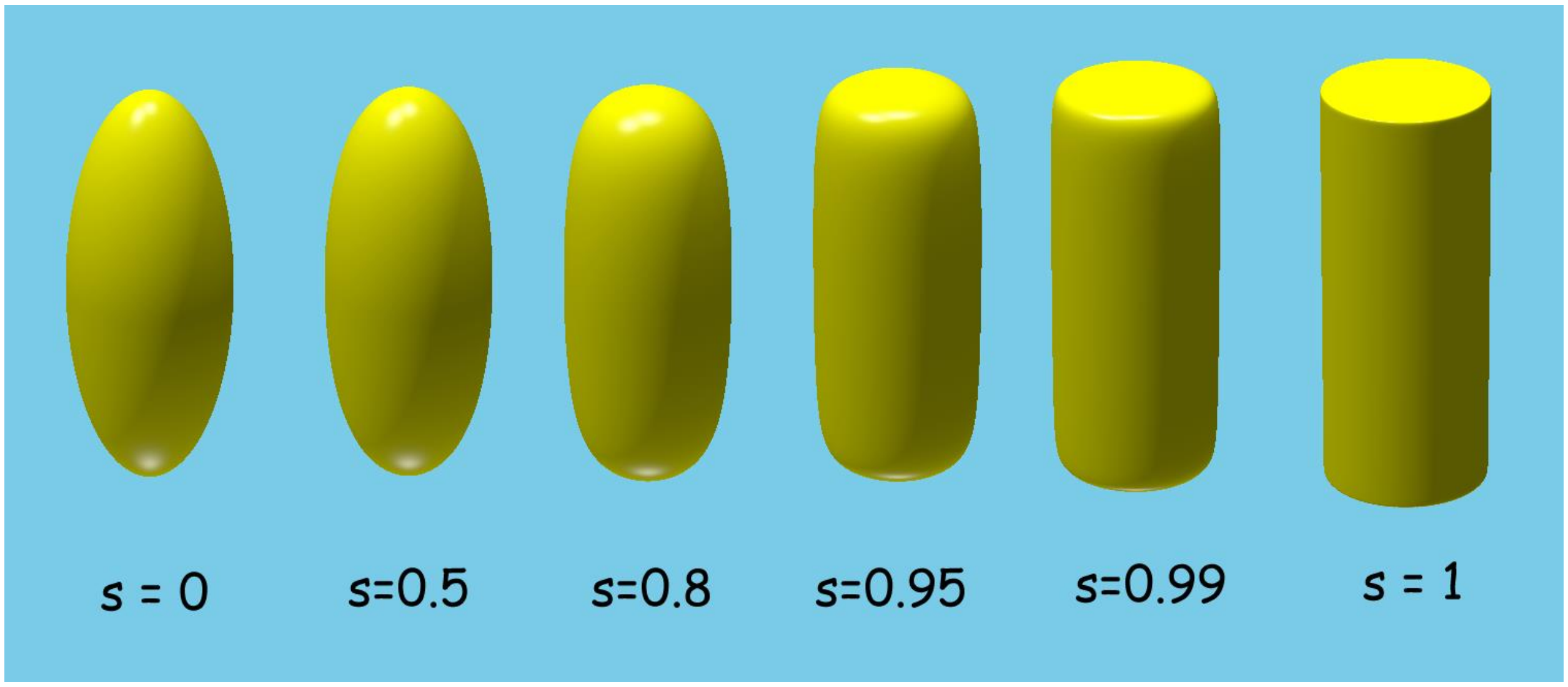

Figure 17: The eccentric sphylinder at varying squareness values
with a=b=1.0 and c=2.5

Just like with the sphube, it is important to subject the sphylinder to these inequalities bounding x, y, and z in order to prune away unwanted lobes.

$$\frac{x^2}{a^2}+\frac{y^2}{b^2}\leq 1$$

$$-c\leq z\leq c$$

The right circular cylinder form of the sphylinder is most useful in practical applications. This occurs when the cylindrical base is a circle with $r=a=b$. For this case, the sphylinder has parametric equations

$$x=r\cos\theta\cos\phi$$

$$y=r\cos\theta\sin\phi$$

$$z=\frac{h\,\sin\theta}{\sqrt{1-s_p\,\cos^2\theta}}$$

where $2h$ is the height of the sphylinder and $s_p$ is a squareness value for the parametric equation.

## 15.8 Cylinder with End Caps

The traditional representation of right circular cylinders in 3D analytic geometry uses the equation

$$x^2 + y^2 = r^2 \qquad with \qquad 0 \le z \le H$$

where $r$ is the radius of the circular base and $H$ is the height of the cylinder. However, this representation has a noticeable limitation in that the end caps of the cylinder are missing. In some applications, it is desirable to include the end caps. Moreover, this equation for the cylinder is disingenuous in three dimensions because the z variable is disjoint from the main equation. In and of itself, this equation is quasi-3D.

Meanwhile, it is possible to represent the cylinder with end caps by using the equation for the sphylinder with $s = 1$. In other words, the cylinder can be represented by the equation

$$\frac{x^2 + y^2}{r^2} + \frac{z^2}{h^2} - \frac{x^2 z^2}{r^2 h^2} - \frac{y^2 z^2}{r^2 h^2} = 1$$

with

$$x^2 + y^2 \le r^2$$
$$-h \le z \le h$$

where $r$ is the radius of the circular base, and $2h$ is the height of the cylinder. This equation provides a genuine 3D representation of the cylinder with end caps. Furthermore, this algebraic equation is a low degree quartic.

## 15.9 Online Resources

The sphube can be viewed interactively using the online graphing calculator Desmos:
https://www.desmos.com/3d/chxxg8y5ih

Likewise, the sphylinder and its cylindrical special case at $s=1$ can be viewed interactively using the online graphing calculator Desmos:
https://www.desmos.com/3d/nopiuowpkg

Some animations of the 3D squircular surfaces discussed in the previous sections were presented at the Joint Mathematics Meetings 2018 conference. The videos are available in:
https://youtu.be/f1nN9DfhgbU

## Summary

In this paper, we provided several formulas related to the FG-squircle. We also presented a three dimensional counterpart of the FG-squircle. The equations are summarized in the following tables.

| Property | Equation(s) |
|---|---|
| FG-Squircle | $x^2 + y^2 - \frac{s^2}{r^2}\,x^2y^2 = r^2$ |
| Polar equation | $\rho = \frac{r\sqrt{2}}{\sqrt{1+\sqrt{1-s^2\sin^2 2\theta}}}$ |
| Area (incomplete type) | $A = 4\frac{r^2}{s}E(\sin^{-1}s, \frac{1}{s})$ |
| Area (complete type) | $A = 4\frac{r^2}{s^2}\,[\,E(s) + (s^2-1)K(s)]$ |
| Closed parametric equation#1 | $x = \frac{r}{2s}\sqrt{2+2s\sqrt{2}\,\cos t + s^2\cos 2t} - \frac{r}{2s}\sqrt{2-2s\sqrt{2}\,\cos t + s^2\cos 2t}$<br>$y = \frac{r}{2s}\sqrt{2+2s\sqrt{2}\,\sin t - s^2\cos 2t} - \frac{r}{2s}\sqrt{2-2s\sqrt{2}\,\sin t - s^2\cos 2t}$ |
| Closed parametric equation#2 | $x = \frac{r\ sgn(\cos t)}{s\sqrt{2}}\sqrt{2+s^2\cos 2t - \sqrt{(2+s^2\cos 2t)^2 - 8\,s^2\cos^2 t}}$<br>$y = \frac{r\ sgn(\sin t)}{s\sqrt{2}}\sqrt{2-s^2\cos 2t - \sqrt{(2-s^2\cos 2t)^2 - 8\,s^2\sin^2 t}}$ |
| Open parametric equation#1 | $x = \frac{r\cos t}{\sqrt{1-s_1\sin^2 t}}$<br>$y = \frac{r\sin t}{\sqrt{1-s_1\cos^2 t}}$<br>$where\ \ s_1 = s\sqrt{2-s^2}$ |
| Open parametric equation#2 | $x = r\cos t$<br>$y = \frac{r\sin t}{\sqrt{1-s^2\cos^2 t}}$ |
| Open parametric equation#3 | $x = \frac{r\cos t}{\sqrt{1-s^2\sin^2 t}}$<br>$y = r\sin t$ |

| 3D Surface | Equation(s) |
|---|---|
| Sphube | $x^2 + y^2 + z^2 - \frac{s^2}{r^2}\, x^2y^2 - \frac{s^2}{r^2}\, y^2z^2 - \frac{s^2}{r^2}\, x^2z^2 + \frac{s^4}{r^4} x^2y^2z^2 = r^2$ |
| Sphube parametric equations | $x = \frac{r \cos\theta \cos\phi}{\sqrt{1 - s_p \cos^2\theta \sin^2\phi - s_p \sin^2\theta}}$ <br> $y = \frac{r \cos\theta \sin\phi}{\sqrt{1 - s_p \cos^2\theta \cos^2\phi - s_p \sin^2\theta}}$ <br> $z = \frac{r \sin\theta}{\sqrt{1 - s_p \cos^2\theta}}$ |
| Squircular ellipsoid | $\frac{x^2}{a^2} + \frac{y^2}{b^2} + \frac{z^2}{c^2} - \frac{s^2x^2y^2}{a^2b^2} - \frac{s^2x^2z^2}{a^2c^2} - \frac{s^2y^2z^2}{b^2c^2} + \frac{s^4x^2y^2z^2}{a^2b^2c^2} = 1$ |
| Squircular cone | $x^2z^2 + y^2z^2 - s^2c^2x^2y^2 - \frac{z^4}{c^2} = 0$ |
| Non-uniform Squircular cone | $\frac{x^2}{a^2} + \frac{y^2}{b^2} - \frac{x^2y^2}{a^2b^2} = \frac{z^2}{c^2}$ |
| Squircular cushion | $x^2 + y^2 + z^2 - \frac{s^2}{r^2}\, x^2y^2 = r^2$ <br> and its rotated variants <br> $x^2 + y^2 + z^2 - \frac{s^2}{r^2}\, y^2z^2 = r^2$ <br> $x^2 + y^2 + z^2 - \frac{s^2}{r^2}\, x^2z^2 = r^2$ |
| Non-uniform Squircular cylinder (no end caps) | $\frac{x^2}{a^2} + \frac{y^2}{b^2} - \frac{z^2x^2y^2}{a^2b^2c^2} = 1$ |

| 3D Surface | Equation(s) |
| --- | --- |
| Sphylinder | $$x^2+y^2+z^2-\frac{s^2}{r^2}\,x^2z^2-\frac{s^2}{r^2}\,y^2z^2=r^2$$ with $$x^2+y^2\le r^2$$ $$-r\le z\le r$$ |
| Rotated variants of the Sphylinder | $$x^2+y^2+z^2-\frac{s^2}{r^2}\,x^2y^2-\frac{s^2}{r^2}\,y^2z^2=r^2$$ $$x^2+y^2+z^2-\frac{s^2}{r^2}\,x^2y^2-\frac{s^2}{r^2}\,x^2z^2=r^2$$ |
| Eccentric sphylinder | $$\frac{x^2}{a^2}+\frac{y^2}{b^2}+\frac{z^2}{c^2}-\frac{s^2}{a^2c^2}\,x^2z^2-\frac{s^2}{b^2c^2}\,y^2z^2=1$$ |
| Stretched sphylinder | $$\frac{x^2+y^2}{r^2}+\frac{z^2}{h^2}-\frac{s^2}{r^2h^2}\,x^2z^2-\frac{s^2}{r^2h^2}\,y^2z^2=1$$ |
| Sphylinder parametric equations | $$x=r\cos\theta\cos\phi$$ $$y=r\cos\theta\sin\phi$$ $$z=\frac{h\ \sin\theta}{\sqrt{1-s_p\ \cos^2\theta}}$$ |
| Cylinder with end caps | $$\frac{x^2+y^2}{r^2}+\frac{z^2}{h^2}-\frac{x^2z^2}{r^2h^2}-\frac{y^2z^2}{r^2h^2}=1$$ with $$x^2+y^2\le r^2$$ $$-h\le z\le h$$ |

## Acknowledgements

We would like to thank Douglas Norton, moderator of the SIGMAA-ARTS, for his tireless efforts and funny jokes during the Joint Mathematics Meetings 2018 in San Diego. We would also like to thank to Bernd Frassek for his suggestions and corrections to this paper. Finally, special thanks to Manuel Fernandez-Guasti for his encouragement and generous praise of this paper.

# Appendix 1

We include Matlab source code for the volumetric mapping in sections 13 and 14.

---

```
function [u,v,w]=sphubeuvw(x,y,z)
% volumetric mapping from cube to sphere
% This is a radial mapping based on the sphube
% input: point (x,y,z) inside cube
% output: point (u,v,w) inside sphere
% Canonical mapping space:
%    (u,v,w) is inside the sphere centered at the origin with radius 1.0
%    (x,y,z)  is inside the square centered at the origin with side length 2.0
x2 = x*x;
y2 = y*y;
z2 = z*z;
r2 = x2 + y2 + z2;
if abs(r2) < eps
  u=x;
  v=y;
  w=z;
  return;
end
s2 = r2 - x2*y2 -y2*z2 - x2*z2 + x2*y2*z2;
m = sqrt(s2)/sqrt(r2);
u = x * m;
v = y * m;
w = z * m;
```

```
function [x,y,z]=sphubexyz(u,v,w)
% volumetric mapping from sphere to cube
% This is a radial mapping based on the sphube
% input: point (u,v,w) inside sphere
% output: point (x,y,z) inside cube
if abs(u) < eps
  x=0;
  [y,z]=fgsquircular(v,w);
  return;
end
if abs(v) < eps
  y=0;
  [x,z]=fgsquircular(u,w);
  return
end
if abs(w) < eps
  z=0;
  [x,y]=fgsquircular(u,v);
  return
end
v2 = v*v;
u2 = u*u;
w2 = w*w;
u4 = u2*u2;
a = v2*w2/u4;
b = -( v2/u2 + w2/u2 + a);
c = 1 + v2/u2 + w2/u2;
d = -(u2 + v2 + w2);
Do = 18*a*b*c*d - 4*b*b*b*d + b*b*c*c - 4*a*c*c*c - 27*a*a*d*d;
if Do<0
  Co = ((2*b*b*b - 9*a*b*c + 27*a*a*d + sqrt(-27*a*a*Do))/2)^(1/3);
  x2 = -1/(3*a) * ( b + Co + (b*b-3*a*c)/Co);
else
  % perform a cube root base on De Moivre's formula
  ro = (2*b*b*b - 9*a*b*c + 27*a*a*d)/2;
  qo = sqrt(27*a*a*Do)/2;
  th=atan2(qo,ro);
  x2 = -1/(3*a) * (b + 2 * (ro*ro+qo*qo)^(1/6) * cos(th/3));
end
x = sign(u)*sqrt(x2);
y = v/u*x;
z = w/u*x;
```

```
function [g,h]=fgsquircular(k,m)
% 2d mapping from circle to square
if abs(k) < eps | abs(m) < eps
  g=k;
  h=m;
  return;
end
sqrt2 = sqrt(2);
u2 = k*k;
v2 = m*m;
uv = k*m;
u2pv2 = u2+v2;
fouru2v2 = 4.*uv.*uv;
rad = u2pv2.*(u2pv2-fouru2v2);
x2 = (u2 + v2 - sqrt(rad))./(2.*v2);
y2 = (u2 + v2 - sqrt(rad))./(2.*u2);
g= sign(uv) * sqrt(x2);
h= sign(uv) * sqrt(y2);
```

## Appendix 2: An Oblique Squircle

In addition to the Lamé curve and the Fernandez-Guasti squircle, there are many other types of squircles. In essence, any shape that can interpolate between the square and the circle can be considered as a squircle. For example, Marc Frantz introduced his own squircle in the October 2018 issue of Mathematics Magazine with the article "*Save This Family! An Instructive Sequence of Curves*".

In this appendix, we will introduce yet another type of squircle, which we will refer to as the *oblique squircle*. The barebones equation for the oblique squircle is

$$\boldsymbol{\cos(\omega x) + \cos(\omega y) = 1 + \cos\omega}$$

This implicit equation becomes a circle as $\omega \to 0$. Similarly, when $\omega = \pi$, the implicit equation produces a square. In between, the implicit equation produces a shape that interpolates between the circle and the square. To be more specific, as $\omega$ goes to zero, the graph of this implicit equation becomes a circle centered at the origin with radius 1. When $\omega$ has value $\pi$, the graph consists of a grid of diagonal lines with a central square at the origin. The square is rotated by 45° and not axis-aligned. This is the reason why the shape is called an oblique squircle

The equation for the oblique squircle can be manipulated to get isolated expressions for both *x* and *y*.

$$x = \frac{\cos^{-1}(1 + \cos\omega - \cos(\omega y))}{\omega} \qquad\qquad y = \frac{\cos^{-1}(1 + \cos\omega - \cos(\omega x))}{\omega}$$

### A2.1 Circular Case when $\omega = 0$

In order to demonstrate that this curve is a type of squircle, we have to first show that as $\omega$ goes to zero, the equation becomes circular. Note that if we just set $\omega = 0$, the equation becomes degenerate of the form *2 = 2*, which cannot be plotted. In order to avoid this degeneracy, we need to take the limit as $\omega$ goes to zero instead. In other words, we need to show that

$$\lim_{\omega \to 0} x \ = \sqrt{1 - y^2}$$

which is equivalent to the equation for the unit circle centered at the origin.

Substituting the expression for x, we get this equation to be proved

$$\lim_{\omega \to 0} \frac{\cos^{-1}(1 + \cos\omega - \cos(\omega y))}{\omega} = \sqrt{1 - y^2}$$

Using L'Hôpital's rule, we can simplify the limit by differentiating both the numerator and denominator with respect to $\omega$ in order to get

$$\lim_{\omega \to 0} \frac{\cos^{-1}(1 + \cos\omega - \cos(\omega y))}{\omega} \ = \ \lim_{\omega \to 0} \frac{\sin\omega - y\ \sin(\omega y)}{\sqrt{1 - (1 + \cos\omega - \cos(\omega y))^2}}$$

Furthermore, we can do a Taylor series expansion of the numerator and denominator at $\omega = 0$ to further simplify.

For the numerator:

$$\sin\omega - y\ \sin(\omega y) = \ \omega(1 - y^2) + \tfrac{1}{6}\omega^3(y^4 - 1) + \tfrac{1}{120}\omega^5(1 - y^6) + O(\omega^7)$$

where $O(\omega^7)$ represents an infinite sum of higher order polynomial terms in $\omega$ with exponents ≥ 7.

For the denominator:

$$1 - (1 + \cos\omega - \cos(\omega y))^2 = \omega^2(1 - y^2) + \tfrac{\omega^4}{6}(-y^4 + 3y^2 - 2) + \tfrac{\omega^6}{360}(14y^6 - 15y^4 - 15y^2 + 16) + O(\omega^8)$$

where $O(\omega^8)$ represents an infinite sum of higher order polynomial terms in $\omega$ with exponents ≥ 8. Substituting in the Taylor series expansions for the numerator and denominator, we get

$$\lim_{\omega\to 0}\frac{\sin\omega - y\ \sin(\omega y)}{\sqrt{1-(1+\cos\omega-\cos(\omega y))^2}}$$

$$= \lim_{\omega\to 0}\frac{\omega(1-y^2)+\frac{1}{6}\omega^3(y^4-1)+\frac{1}{120}\omega^5(1-y^6)+O(\omega^7)}{\sqrt{\omega^2(1-y^2)+\frac{1}{6}\omega^4(-y^4+3y^2-2)+\frac{1}{360}\omega^6(14y^6-15y^4-15y^2+16)+O(\omega^8)}}$$

$$= \lim_{\omega\to 0}\frac{\omega(1-y^2)\quad\left(1+\frac{1}{6}\frac{\omega^3(y^4-1)}{\omega(1-y^2)}+\frac{1}{120}\frac{\omega^5(1-y^6)}{\omega(1-y^2)}+\frac{1}{\omega(1-y^2)}O(\omega^7)\right)}{\sqrt{\omega^2(1-y^2)}\sqrt{1+\frac{1}{6}\frac{\omega^4(-y^4+3y^2-2)}{\omega^2(1-y^2)}+\frac{1}{360}\frac{\omega^6(14y^6-15y^4-15y^2+16)}{\omega^2(1-y^2)}+\frac{1}{\omega^2(1-y^2)}O(\omega^8)}}$$

$$= \lim_{\omega\to 0}\frac{\sqrt{1-y^2}\quad\left(1-\frac{1}{6}\omega^2(1+y^2)+\frac{1}{120}\omega^4(1+y^2+y^4)+\frac{1}{(1-y^2)}O(\omega^6)\right)}{\sqrt{1+\frac{1}{6}\omega^2(y^2-2)+\frac{1}{360}\omega^4(-14y^4+y^2+16)+\frac{1}{(1-y^2)}O(\omega^6)}}$$

$$= \sqrt{1-y^2}$$

### A2.2 Square Case when $\omega = \pi$

We now consider the other limiting case at $\omega = \pi$. The oblique squircle equation reduces to:

$$\cos(\pi x)+\cos(\pi y)=0$$

This means

$$\cos(\pi y) = -\cos(\pi x)$$

$$\Rightarrow \qquad \cos(\pi y)=\cos(\pi x+(2n+1)\pi)$$

where n is any integer.

This simplifies to

$$\pi y = \pm\pi x+(2n+1)\pi$$

$$\Rightarrow \qquad y=\pm x+2n+1$$

This means that the graph consists of many slanted lines with a slope of ±1. At the center of these multitude of lines lays a tilted square.

### A2.3 Equation with Scale and Squareness

It is convenient to incorporate scale into the equation of the oblique squircle. Thus, we have this alternative equation generalized from the barebones equation previously provided. The variable *s* represents the squareness of the oblique squircle. When *s=1*, the shape is a square with side length *2k*. As *s ⟶0*, the shape becomes a circle with radius *k*.

$$\boldsymbol{\cos\left(\frac{s\pi x}{k}\right)+\cos\left(\frac{s\pi y}{k}\right)=1+\cos(s\pi)}$$

Note that the graph of the oblique squircle actually extends out to span the whole 2D plane. In fact, the oblique squircle is doubly-periodic in the plane. This is because both *x* & *y* can have values in (-∞, +∞) and the cosine function is a periodic function. However, our main region of interest is confined to $-k \le x \le k$ and $-k \le y \le k$.

### A2.4 Equation with Offsets and Eccentricities

The equation for the oblique squircle can be generalized further to incorporate center offsets and eccentricities. Suppose we want our shape to centered at the point $(x_0, y_0)$ instead of the origin. Furthermore, suppose we want anisotropic scale factors of $k_x$ and $k_y$. The generalized equation for the oblique squircle becomes

$$\cos\left(\frac{s\pi(x - x_0)}{k_x}\right) + \cos\left(\frac{s\pi(y - y_0)}{k_y}\right) = 1 + \cos(s\pi)$$

When $s=1$, the central shape is a rhombus. As $s \to 0$, the shape becomes an ellipse with semi-major and semi-minor axial lengths of $k_x$ and $k_y$.

### A2.5 Three-Dimensional Counterpart

Just like the Fernandez-Guasti squircle, the oblique squircle can be extended to 3D. However, the corresponding 3D shape to the square will not be a cube. Instead, it will be a shape that resembles the octahedron which we shall name as the *sham octahedron*. The equation for the 3D implicit surface is

$$\cos\left(\frac{s\pi x}{k}\right) + \cos\left(\frac{s\pi y}{k}\right) + \cos\left(\frac{s\pi z}{k}\right) = 2 + \cos(s\pi)$$

This implicit equation becomes a sphere as $s \to 0$. Similarly, when $s = 1$, the implicit equation produces the sham octahedron. In between, the implicit equation produces a octahedral shape that interpolates between the sphere and the sham octahedron. The other parameter $k$ provides a scale factor in the shapes.

### Addendum

For more information on the oblique squircle and related shapes, see the follow-up paper by the same author, titled "*Visualizing Squircular Implicit Surfaces*". It was presented at the Joint Mathematics Meetings in 2024.

https://arxiv.org/abs/2210.15232